\documentclass[11pt,british]{amsart}

\newif\ifpdf
\ifx\pdfoutput\undefined\pdffalse\else\pdfoutput=1\pdftrue\fi\newif\pdf
\ifpdf\relax\else
\usepackage[T1]{fontenc}
\fi

\usepackage{babel,eucal,url,amssymb,stmaryrd,
enumerate,amscd,
}

\usepackage[pagebackref]{hyperref}

\theoremstyle{plain}
\newtheorem{thm}{Theorem}[section]
\newtheorem{lem}[thm]{Lemma}
\newtheorem{pro}[thm]{Proposition}
\newtheorem{co}[thm]{Corollary}

\theoremstyle{definition}

\theoremstyle{remark}
\newtheorem{rem}[thm]{Remark}

\newcommand{\Gtwo}{\ifmmode{{\rm G}_2}\else{${\rm G}_2$}\fi}

\date{\today}
\begin{document}

\title[Geometry of Paraquaternionic K\"ahler manifolds with torsion]%
{Geometry of Paraquaternionic K\"ahler manifolds with torsion}

\author{Simeon Zamkovoy}
\address[Zamkovoy]{University of Sofia "St. Kl. Ohridski"\\
Faculty of Mathematics and Informatics,\\
Blvd. James Bourchier 5,\\
1164 Sofia, Bulgaria} \email{zamkovoy@fmi.uni-sofia.bg}

\begin{abstract}
We study the geometry of PQKT-connections. We find conditions to
the existence of a PQKT-connection and prove that if it exists it
is unique. We show that PQKT geometry persist in a conformal class
of metrics. \\[8mm]

MSC: 53C15, 5350, 53C25, 53C26, 53B30
\end{abstract}

\maketitle \setcounter{tocdepth}{2} \tableofcontents

\section{Introduction and statement of the results}

 The geometry of locally supersymmetric vector multiplets in
dimension $(1,3)$ is known as projective special K\"ahler
geometry. As indicated by the name, such manifolds can be obtained
from affine special K\"ahler manifolds (with homogeneity
properties) by a projectivization. This can also be understood
from the physical point of view in terms of the conformal
calculus. One first constructs a superconformally invariant theory
and then eliminates conformal compensators by imposing gauge
conditions. This gauge fixing amounts to the projectivization of
the scalar manifold underlying the superconformal theory
\cite{W1,S1,C1,C2,A3,A4,C3,C6,F3,A5,C4}. Based upon this results
such a construction could be adapted to the case of Euclidean
signature and projective special paraK\"ahler manifolds are
constructed in \cite{C4}. In Minkowski signature the coupling to
supergravity implies that the scalar geometry is
quaternionic-K\"ahler instead of hyper-K\"ahler manifolds
\cite{B5}. The relation between these two kinds of geometries can
again be understood as projectivization, because every
quaternionic-K\"ahler manifold can be obtained as the quotient of
a hyper-K\"ahler cone \cite{S2}. The identification of the scalar
geometry of Euclidean hypermultiplets in rigidly supersymmetric
theories is made in \cite{C5}. The fact, that the scalar manifolds
of Euclidean hypermultiplets are hyper-paraK\"ahler manifolds is
one of the main results in \cite{C5}.

 We recall that an almost hyper paracomplex structure on a 4n-dimensional manifold
$M$ is a triple $H=(J_{\alpha}), \alpha=1,2,3$, of two almost
paracomplex structures and one complex structure
$J_{\alpha}:TM\rightarrow TM$ satisfying the paraquaternionic
identities
$$J_{\alpha}^2=\epsilon_{\alpha},J_{\alpha}J_{\beta}=-J_{\beta}J_{\alpha}=-\epsilon_{\gamma}J_{\gamma},\alpha,\beta,\gamma
= 1,2,3, \epsilon_1=\epsilon_2=-\epsilon_3=1.$$
 Here and henceforth $(\alpha,\beta,\gamma)$ is a cyclic
permutation of $(1,2,3)$.

 When each $J_{\alpha}$ is an integrable almost (para)complex structure,
$H$ is said to be a hyper paracomplex structure on $M$. Such a
structure is also called sometimes pseudo-hyper-complex \cite{D1}.
Any hyper paracomplex structure admits a unique torsion-free
connection $\nabla^{CP}$ preserving $J_1,J_1,J_3$ \cite{A1,A2}
called the complex product connection. Examples of
hyper-paracomplex structures on the simple Lie groups
$SL(2n+1,\mathbb R)$,$SU(n,n+1)$ are constructed in \cite{I3}.

 Almost paraquaternionic structures are introduced by Libermann
(Libermann called them almost quaternionic structure of the second
kind) \cite{L1}. An almost paraquaternionic structure on $M$ is a
rank-3 subbundle $\mathbb P$ $\subset End(TM)$ which is locally
spanned by an almost hyper-paracomplex structure $H=(J_{\alpha})$;
such a locally defined triple $H$ will be called an admissible
basis of $\mathbb P$. A linear connection $\nabla$ on $TM$ is
called paraquaternionic if $\nabla$ preserves $\mathbb P$, i.e.
$\nabla_X\sigma\in \Gamma (\mathbb P)$ for all vector fields $X$
and smooth sections $\sigma \in \Gamma (\mathbb P)$. An almost
paraquaternionic structure is said to be a paraquaternionic if
there is a torsion-free paraquaternionic connection. A $\mathbb
P$-hermitian metric is a pseudo-Riemannian metric which is
compatible with $\mathbb P$ in the sense that
$g(J_{\alpha}X,J_{\alpha}Y)=-\epsilon_{\alpha}g(X,Y)$, for
$\alpha=1,2,3$. The signature of $g$ is necessarily of neutral
type $(2n,2n)$. An almost paraquaternionic (resp.
paraquaternionic) manifold with a $\mathbb P$-hermitian metric is
said to be an almost paraquaternionic Hermitian (resp.
paraquaternionic Hermitian) manifold.

For $n\geq 2$, the existence of a torsion-free paraquaternionic
connection is a strong condition which is equivalent to the
1-integrability of the associated $GL(n,H)Sp(1,\mathbb R)$
-structure \cite{A1,A2}. The paraquaternionic condition controls
the Nijenhuis tensor in the sense that
$N(X,Y)J_{\alpha}:=N_{\alpha}(X,Y)$ preserves the subbundle
$\mathbb P$. An invariant first order differential operator $D$ is
defined on any almost paraquaternionic manifolds which is two-step
nilpotent i.e. $D^2=0$ exactly when the structure is
paraquaternionic \cite{I4}. If the Levi-Civita connection of a
paraquaternionic hermitian manifold $(M,g,\mathbb P)$ is a
paraquaternionic connection then $(M,g,\mathbb P)$ is called
Paraquaternionic K\"ahler (briefly {PQK}). This condition is
equivalent the statement that the holonomy group of $g$ is
contained in $Sp(n,\mathbb R)$.$Sp(1,\mathbb R)$ \cite{I3}. A
typical example is the paraquaternionic projective space endowed
with the standard paraquaternionic K\"ahler structure \cite{B1}.
Any paraquaternionic K\"ahler manifold of dimension $4n\geq 8$ is
known to be Einstein \cite{G1,V1}. If a PQK manifold there exist a
global admissible basis $(H)$ such that each almost (para)complex
structure $(J_{\alpha})\in(H), \alpha =1,2,3$ is parallel with
respect to the Levi-Civita connection then the manifold is called
hyper para-K\"ahler (briefly HPK). Such manifolds are also called
hypersymplectic \cite{H1}, neutral hyper-K\"ahler \cite{F1,K1}. In
this case the holonomy group of $g$ is contained in $Sp(n,\mathbb
R)$, $n>1$ \cite{V1}. Twistor and reflector spaces on
paraquaternionic K\"ahler manifold are constructed and the
integrability of the associated (para)complex structures are
investigated in \cite{B4} and \cite{J1}, respectively. These
constructions work also in the paraquaternionic case \cite{I5}.

A natural generalization of PQK spaces is the notion of
paraquaternionic K\"ahler manifolds with torsion (briefly PQKT)
which means that there exists a paraquaternionic connection
preserving the metric g with totally skew-symmetric torsion of
type $(1,2)+(2,1)$ with respect to each $J_{\alpha}$. More
general, if one considers the same construction in the general
case of an almost paraquaternionic structure and defines the
almost complex structure on the twistor space (resp. almost
paracomplex structure on the reflector space) using horizontal
spaces of an arbitrary linear connection than the integrability
condition is equivalent to the condition that the torsion is of
type $(0,2)$ with respect each $J_{\alpha}$ \cite{I6}. The main
object of interest in this article is the differential geometric
properties of PQKT manifolds.

In section 2 we find necessary and sufficient conditions for the
almost paraquaternionic structure to be paraquaternionic structure
if the dimension is at least 8 (Theorem~\ref{thm0.4}).

In section 3 we find necessary and sufficient conditions to the
existence of a PQKT connection in terms of the K\"ahler 2-forms
and  show that the PQKT-connection is unique if the dimension is
at least 8 (Theorem~\ref{thm1}) and we prove that the PQKT
manifolds are invariant under conformal transformations of the
metric.

In section 4 we prove that the $(2,0)+(0,2)$-parts of the Ricci
forms $\rho_{\alpha}$, $\rho_{\beta}$ with respect to $J_{\gamma}$
coincide if the dimension is at least 8 (Theorem~\ref{thm10}). We
define torsion 1-form t as a suitable trace of the torsion 3-form
and show that $\rho_{\alpha}$ is of type $(1,1)$ with respect to
$J_{\alpha}$ if and only if $dt$ is of type $(1,1)$ with respect
to each $J_{\alpha}$,$\alpha=1,2,3$ provided the dimension is at
least 8 (Theorem~\ref{thm10.1}). We show that $\star$-Ricci tensor
$\rho^{\star}_{\alpha}$ is symmetric if and only if $dt$ is of
type $(1,1)$ with respect to each $J_{\alpha}$,$\alpha=1,2,3$
(Corollary~\ref{10.2}).

In section 5 we show that there are no homogeneous proper PQKT
manifolds (i.e. homogeneous PQKT which is not PQK or HPKT) with
$dT$ of type (2,2) provided that the torsion is parallel and
dimension is at least 8 (Theorem~\ref{thm2}).

\section{Almost Paraquaternionic structures}

 In this section we study of the integrability of
almost paraquaternionic structures. We recall that an almost
paraquaternionic structure is G-structure with structure group
$GL(n,H)$.$Sp(1,\mathbb R)$$\cong $\\$GL(2n,\mathbb
R)$.$Sp(1,\mathbb R)$. The almost paraquaternionic structures are
studied intensively in the last years, especially in the case of
hyper-paraK\"ahler and paraquaternionic K\"ahler
\cite{G1,V1,B1,B2,B3,I5} i.e. the structure group is further
reduced to $Sp(n,\mathbb R)$ or $Sp(n,\mathbb R)$.$Sp(1,\mathbb
R)$.

 Let $\nabla$ be a paraquaternionic connection i.e.
\begin{equation}\label{1}
\nabla J_{\alpha}=\omega_{\beta}\otimes J_{\gamma} +
\epsilon_{\gamma} \omega_{\gamma}\otimes J_{\beta},
\end{equation}
where $\omega_{\alpha}, \alpha =1,2,3$ are 1-forms.

 Here and henceforth $(\alpha, \beta, \gamma)$ is a cyclic
permutation of $(1,2,3)$.

 The Nijenhuis tensor $N_{\alpha}$ of an almost (para)complex
structure $J_{\alpha}$ is given by
\\
$N_{\alpha}(X,Y)=[J_{\alpha}X,J_{\alpha}Y]+\epsilon_{\alpha}[X,Y]-J_{\alpha}[J_{\alpha}X,Y]
- J_{\alpha}[X,J_{\alpha}Y].$

 It is well-known that an almost (para)complex structure is a (para)complex
structure if and only if its Nijenhuis tensor vanishes.

 The Nijenhuis bracket $[[A,B]]$ of two endomorphisms is defined in
terms of the Lie bracket of vector fields in the following way:
$$[[A,B]](X,Y)=[AX,BY]-A[BX,Y]-B[X,AY]+[BX,AY]-B[AX,Y]-A[X,BY]+$$
$$+(AB+BA)[X,Y]$$
and $[[J_{\alpha},J_{\alpha}]](X,Y)=2N_{\alpha}(X,Y)$.

\begin{pro}\label{thm0.1}
Let $H=(J_{\alpha})$ be an almost hyper-paracomplex structure on
$M$.

1) There exist connection $\nabla^{CP}$ which preserves $H$. The
connection $\nabla^{CP}$ is given by
\begin{equation}\label{tt1.1}
\nabla^{CP}_{X}Y=\frac{1}{12}(\sum_{(\alpha,\beta,\gamma)}(J_{\alpha}[J_{\beta}X,J_{\gamma}Y]-J_{\alpha}[J_{\gamma}X,J_{\beta}Y])-2\sum_{\alpha=1}^{3}(\epsilon_{\alpha}J_{\alpha}[J_{\alpha}X,Y]-\epsilon_{\alpha}J_{\alpha}[X,J_{\alpha}Y]))-
\end{equation}
$$-\frac{1}{12}\sum_{\alpha=1}^{3}(\epsilon_{\alpha}[J_{\alpha}X,J_{\alpha}Y]-\epsilon_{\alpha}J_{\alpha}[J_{\alpha}X,Y]-\epsilon_{\alpha}J_{\alpha}[X,J_{\alpha}Y]+[X,Y])+\frac{1}{2}[X,Y],$$
where $(\alpha,\beta,\gamma)$ indicates sum over cyclic
permutations of $(1,2,3)$.

2) its torsion tensor $T^H$ satisfies
\begin{equation}\label{tt1.2}
T^H=-\frac{1}{12}\sum_{\alpha=1}^{3}\epsilon_{\alpha}[[J_{\alpha},J_{\alpha}]],
\end{equation}
for any two vector fields $X,Y$.
\end{pro}
{\it Proof.} One verifies easily that the formula $(\ref{tt1.1})$
defines a connection $\nabla^{CP}$ which preserve $H$ and whose
torsion tensor is given by $(\ref{tt1.2})$. \hfill {\bf Q.E.D.}

\begin{rem}
If we denote by
$$\nabla^{0}_{X}Y=\frac{1}{12}(\sum_{(\alpha,\beta,\gamma)}(J_{\alpha}[J_{\beta}X,J_{\gamma}Y]-J_{\alpha}[J_{\gamma}X,J_{\beta}Y])-2\sum_{\alpha=1}^{3}(\epsilon_{\alpha}J_{\alpha}[J_{\alpha}X,Y]-\epsilon_{\alpha}J_{\alpha}[X,J_{\alpha}Y]))+\frac{1}{2}[X,Y],$$
then the connection $\nabla^0$ is characterized by the following
properties:

1) It is torsion-free: $T^{\nabla^0}=0$;

2) $\nabla^0J_{\alpha}=-\frac{1}{2}[T^H,J_{\alpha}]$.

The $\nabla^{CP}$ connection can be written as
$$\nabla^{CP}=\nabla^{0}+\frac{1}{2}T^H.$$
\end{rem}
 We need from the following
\begin{lem}
Let $H=(J_{\alpha})$ be an almost hyper-paracomplex structure on
$M$. Then for any vector $X,Y$ the following formulas hold:
\begin{equation}\label{l0.2}
[[J_{\alpha},J_{\beta}]](X,Y)=J_{\alpha}T^H(X,J_{\beta}Y)+J_{\alpha}T^H(J_{\beta}X,Y)+J_{\beta}T^H(X,J_{\alpha}Y)+J_{\beta}T^H(J_{\alpha}X,Y)-
\end{equation}
$$+J_{\beta}T^H(J_{\alpha}X,Y)-T^H(J_{\alpha}X,J_{\beta}Y)-T^H(J_{\beta}X,J_{\alpha}Y);$$
\begin{equation}\label{l0.3}
-12\epsilon_{\gamma}T^H(X,Y)=J_{\alpha}[[J_{\alpha},J_{\beta}]](X,J_{\beta}Y)+J_{\alpha}[[J_{\alpha},J_{\beta}]](J_{\beta}X,Y)+J_{\beta}[[J_{\alpha},J_{\beta}]](X,J_{\alpha}Y)+
\end{equation}
$$+J_{\beta}[[J_{\alpha},J_{\beta}]](J_{\alpha}X,Y)-[[J_{\alpha},J_{\beta}]](J_{\alpha}X,J_{\beta}Y)-[[J_{\alpha},J_{\beta}]](J_{\beta}X,J_{\alpha}Y);$$
\begin{equation}\label{l0.4}
\frac{1}{2}[[J_{\alpha},J_{\alpha}]](X,Y)=-\epsilon_{\alpha}T^H(X,Y)+J_{\alpha}T^H(X,J_{\alpha}Y)+J_{\alpha}T^H(J_{\alpha}X,Y)-T^H(J_{\alpha}X,J_{\alpha}Y);
\end{equation}
\begin{equation}\label{l0.5}
2[[J_{\alpha},J_{\alpha}]](X,Y)=[[J_{\beta},J_{\beta}]](J_{\gamma}X,J_{\gamma}Y)-J_{\gamma}[[J_{\beta},J_{\beta}]](J_{\gamma}X,Y)-J_{\gamma}[[J_{\beta},J_{\beta}]](X,J_{\gamma}Y)-
\end{equation}
$$\epsilon_{\gamma}[[J_{\beta},J_{\beta}]](X,Y)+[[J_{\gamma},J_{\gamma}]](J_{\beta}X,J_{\beta}Y)-J_{\beta}[[J_{\gamma},J_{\gamma}]](J_{\beta}X,Y)-$$
$$-J_{\beta}[[J_{\gamma},J_{\gamma}]](X,J_{\beta}Y)-\epsilon_{\beta}[[J_{\gamma},J_{\gamma}]](X,Y)$$
\end{lem}
{\it Proof.} The first three of these equalities follow by
definition with long but standard computation. The fourth equality
is a Proposition 6.1 in \cite{I5}. \hfill {\bf Q.E.D.}

 As an application of these formulas we obtain necessary and sufficient conditions an
almost hyper-paracomplex structure $H$ to be hyper-paracomplex
structure.
\begin{pro}\label{thm0.2}
Let $H=(J_{\alpha})$ be an almost hyper-paracomplex structure on
$M$. Then the following conditions are equivalent:

1) $H$ is a hyper-paracomplex structure;

2) two of the almost (para)complex structures $J_{\alpha}$
$(\alpha=1,2,3)$ are integrable;

3) one of the Nijenhuis brackets $[[J_{\alpha},J_{\beta}]]$
$(\alpha\not= \beta)$ is zero.

If one of this condition is verified all Nijenhuis brackets
$[[J_{\alpha},J_{\beta}]]$, $\forall \alpha, \beta$ vanish.
\end{pro}
{\it Proof.} If $H$ is hyper-paracomplex then $T^H=0$ and 2), 3)
follow by $(\ref{l0.4})$, $(\ref{l0.2})$ respectively; Vice-versa,
2) or 3) imply 1) by $(\ref{l0.5})$, $(\ref{l0.3})$ respectively.
\hfill {\bf Q.E.D.}
\begin{pro}\label{thm0.3}
Let $\mathbb P$ be an almost paraquaternionic structure and
$H=(J_{\alpha})$ be an admissible basis of $\mathbb P$. Let
$\nabla$ be a globally defined connection which preserve $\mathbb
P$ and let $T$ be its torsion tensor. Then

1) There exist globally defined connection $\nabla^{\mathbb P}$
which preserves $\mathbb P$. The connection $\nabla^{\mathbb P}$
is given by
\begin{equation}\label{tt1.5}
\nabla^{\mathbb
P}_{X}Y=\nabla_XY+\sum_{\alpha=1}^{3}(\epsilon_{\alpha}b_{\alpha}-\frac{1}{3}\epsilon_{\alpha}b\circ
J_{\alpha})(X)J_{\alpha}Y-
\end{equation}
$$-\frac{1}{12}\sum_{\alpha=1}^{3}(T(X,Y)-\epsilon_{\alpha}T(J_{\alpha}X,J_{\alpha}Y)+\epsilon_{\alpha}J_{\alpha}T(X,J_{\alpha}Y)+\epsilon_{\alpha}J_{\alpha}T(J_{\alpha}X,Y)),$$
where $b_{\alpha}$, $b$ are local 1-forms defined by
$$b_{\alpha}(X)=\frac{1}{2n-1}tr(J_{\alpha}T(X))=-\frac{1}{2n-1}\sum_{i=1}^{4n}\epsilon_{i}g(T(X,e_i),J_{\alpha}e_i), b=\sum_{\alpha=1}^{3}\epsilon_{\alpha} b_{\alpha}\circ J_{\alpha},\alpha=1,2,3;$$

2) The torsion tensor $T^{\mathbb P}$ of $\nabla^{\mathbb P}$ is
given by
\begin{equation}\label{tt1.6}
T^{\mathbb P}=T^H + \partial(C^H),
\end{equation}
where
$C^H=\sum_{\alpha=1}^{3}\epsilon_{\alpha}a^{H}_{\alpha}\otimes
J_{\alpha}$, $\partial$ denotes the operator of alternation and
$$a^{H}_{\alpha}(X)=\frac{1}{2n-1}tr(J_{\alpha}T^H(X))=-\frac{1}{2n-1}\sum_{i=1}^{4n}\epsilon_{i}g(T^H(X,e_i),J_{\alpha}e_i),\alpha=1,2,3$$
are the structure 1-forms of $H$.

Moreover
\begin{equation}\label{tt1.3}
\sum_{\alpha=1}^{3}\epsilon_{\alpha}a^{H}_{\alpha}\circ
J_{\alpha}=0.
\end{equation}
\end{pro}
{\it Proof.} For any connection $\nabla$ with torsion tensor $T$
preserving $\mathbb P$ we have
$$[[J_{\alpha},J_{\alpha}]](X,Y)=2\partial ((\epsilon_{\gamma}\omega_{\gamma}\circ J_{\alpha}-\epsilon_{\beta}\omega_{\beta})\otimes J_{\beta}+(\omega_{\beta}\circ J_{\alpha}+\omega_{\gamma})\otimes J_{\gamma})(X,Y)-$$
$$-2T(J_{\alpha}X,J_{\alpha}Y)+2J_{\alpha}T(J_{\alpha}X,Y)+2J_{\alpha}T(X,J_{\alpha}Y)-2\epsilon_{\alpha}T(X,Y).$$
Thus, we have
\begin{equation}\label{tt1.7}
6T^H(X,Y)+\partial \sum_{(\alpha, \beta,
\gamma)}((2\epsilon_{\beta}\omega_{\alpha}+\epsilon_{\beta}\omega_{\gamma}\circ
J_{\beta}-\epsilon_{\alpha}\omega_{\beta}\circ J_{\gamma})\otimes
J_{\alpha})(X,Y)=
\end{equation}
$$=\sum_{\alpha=1}^{3}(T(X,Y)+\epsilon_{\alpha}T(J_{\alpha}X,J_{\alpha}Y)-\epsilon_{\alpha}J_{\alpha}T(J_{\alpha}X,Y)-\epsilon_{\alpha}J_{\alpha}T(X,J_{\alpha}Y)).$$
One verifies easily that the formula $(\ref{tt1.5})$ defines a
connection $\nabla^{\mathbb P}$ which preserve $\mathbb P$ and
whose torsion tensor is
\begin{equation}\label{tt1.8}
6T^{\mathbb
P}(X,Y)=\sum_{\alpha=1}^{3}(T(X,Y)+\epsilon_{\alpha}T(J_{\alpha}X,J_{\alpha}Y)-\epsilon_{\alpha}J_{\alpha}T(J_{\alpha}X,Y)-\epsilon_{\alpha}J_{\alpha}T(X,J_{\alpha}Y))+
\end{equation}
$$+6\partial\sum_{\alpha=1}^{3}((\epsilon_{\alpha}b_{\alpha}-\frac{1}{3}b\circ
J_{\alpha})\otimes J_{\alpha})(X,Y)$$ Taking the appropriate trace
in $(\ref{tt1.7})$, to get
\begin{equation}\label{tt1.9}
3\epsilon_{\alpha}a_{\alpha}^{H}=(2\epsilon_{\beta}\omega_{\alpha}+\epsilon_{\beta}\omega_{\gamma}\circ
J_{\beta}-\epsilon_{\alpha}\omega_{\beta}\circ
J_{\gamma})+(2\epsilon_{\alpha}b_{\alpha}+b_{\beta}\circ
J_{\gamma}-b_{\gamma}\circ J_{\beta}).
\end{equation}
Now, equalities $(\ref{tt1.7})$, $(\ref{tt1.8})$ and
$(\ref{tt1.9})$ prove $(\ref{tt1.6})$. The equality
$(\ref{tt1.3})$ follows from $(\ref{tt1.9})$. \hfill {\bf Q.E.D.}

 The following theorem gives the necessary and
sufficient condition for an almost paraquaternionic structure to
be paraquaternionic structure.
\begin{thm}\label{thm0.4}
An almost paraquaternionic structure $\mathbb P$ is
paraquaternionic structure if and only if $T^{\mathbb P}=0$.
\end{thm}
{\it Proof.} 1) Assume $T^{\mathbb P}=0$. Then $T^H$ has the form
$$T^H=-\partial \sum_{\alpha=1}^{3}\epsilon_{\alpha}a^{H}_{\alpha}\otimes J_{\alpha}.$$
It is easy check that $\nabla=\nabla^H +
\sum_{\alpha=1}^{3}\epsilon_{\alpha}a^{H}_{\alpha}\otimes
J_{\alpha}$ is torsion-free connection. From equalities $\nabla
J_{\alpha}=\bar {\omega_{\beta}}\otimes J_{\gamma}+
\epsilon_{\gamma} \bar {\omega_{\gamma}}\otimes J_{\beta},$ where
$\bar
{\omega_{\alpha}}=\omega_{\alpha}-2\epsilon_{\gamma}a^{H}_{\alpha}$
it follows that $\nabla$ preserves $\mathbb P$.

2) How let $\nabla$ be a paraquaternionic connection and
$H=(J_{\alpha})$ an admissible basis of $\mathbb P$. For any
torsion-free connection $\bar {\nabla}$ preserving $\mathbb P$ we
have
$$[[J_{\alpha},J_{\alpha}]]=2\partial ((\epsilon_{\gamma}\omega_{\gamma}\circ J_{\alpha}-\epsilon_{\beta}\omega_{\beta})\otimes J_{\beta}+(\omega_{\beta}\circ J_{\alpha}+\omega_{\gamma})\otimes J_{\gamma})$$
Thus we obtain
\begin{equation}\label{tt1.10}
6T^H=-\partial \sum_{(\alpha,
\beta,\gamma)}(2\epsilon_{\beta}\omega_{\alpha}+\epsilon_{\beta}\omega_{\gamma}\circ
J_{\beta}-\epsilon_{\alpha}\omega_{\beta}\circ J_{\gamma})\otimes
J_{\alpha}.
\end{equation}
From the formula $(\ref{tt1.9})$, we get $6T^H=-6\partial
\sum_{\alpha=1}^{3}\epsilon_{\alpha}a^{H}_{\alpha}\otimes
J_{\alpha}$. Hence $T^{\mathbb P}=0$. \hfill {\bf Q.E.D.}

\section{Characterizations of PQKT connection}

Let $(M,g,(J_{\alpha})\in \mathbb P$,$\alpha=1,2,3)$ be a
4n-dimensional almost paraquaternionic manifold with $\mathbb
P$-hermitian pseudo-Riemannian metric $g$. We shell work locally
with an admissible basis $(J_{\alpha})$. The K\"ahler form
$F_{\alpha}$ of each $J_{\alpha}$ is defined by
$F_{\alpha}=g(.,J_{\alpha}.)$. The corresponding Lee forms are
given by $\theta_{\alpha}=-\epsilon_{\alpha}\delta F_{\alpha}\circ
J_{\alpha}$.

For an $r$-form $\psi$ we denote by $J_{\alpha}\psi$ the $r$-form
defined by
$J_{\alpha}\psi(X_1,...,X_r):=(-1)^r\psi(J_{\alpha}X_1,...,J_{\alpha}X_r),\\
\alpha=1,2,3$. We shall use the notations
$d_{\alpha}F_{\beta}(X,Y,Z)=-dF_{\beta}(J_{\alpha}X,J_{\alpha}Y,J_{\alpha}Z),
\alpha,\beta =1,2,3.$

We recall the decomposition of a skew-symmetric tensor $P\in
\Lambda^2T^*M \otimes TM$ with respect to a given almost
(para)complex structure $J_{\alpha}$. The (1,1), (2,0) and (0,2)
part of $P$ are defined by $P^{1,1}(J_{\alpha}X,J_{\alpha}Y) = -
\epsilon_{\alpha}P^{1,1}(X,Y), P^{2,0}(J_{\alpha}X,Y) =
J_{\alpha}P^{2,0}(X,Y),
 P^{0,2}(J_{\alpha}X,Y) = - J_{\alpha}P^{0,2}(X,Y)$, respectively.

For each $\alpha = 1,2,3$, we denote by $dF^+_{\alpha}$ (resp.
$dF^-_{\alpha}$) the $(1,2)+(2,1)$-part (resp. $(3,0)+(0,3)$-part)
of  $dF_{\alpha}$ with respect to the almost (para)complex
structure $J_{\alpha}$. We consider the following 1-forms
$$
\theta_{\alpha,\beta} = \epsilon_{\alpha}
\frac{1}{2}\sum_{i=1}^{4n}
\epsilon_{i}dF^+_{\alpha}(X,e_i,J_{\beta}e_i), \quad \alpha ,\beta
=1,2,3.
$$
 Here and further
$e_1,e_2,\dots,e_{n},e_{n+1}=J_3e_1,e_{n+2}=J_3e_2,\dots,e_{2n}=J_3e_n,e_{2n+1}=J_1e_1,e_{2n+2}=J_1e_2,\dots,e_{3n}=J_1e_n,e_{3n+1}=J_2e_1,e_{3n+2}=J_2e_2,\dots,e_{4n}=J_2e_n$
where
$J_1e_i=e_{2n+i},J_1e_{n+i}=e_{3n+i},J_2e_i=e_{3n+i},J_2e_{n+i}=-e_{2n+i},J_3e_{2n+i}=-e_{3n+i},J_3e_i=e_{2n+i},i=1\dots
n$ and $g(e_i,e_i)=\epsilon_i$ where
$\epsilon_i=+1,$$i=1,\dots,2n$ and
$\epsilon_i=-1,$$i=2n+1,\dots,4n$ is an orthonormal basis of the
tangential space.

Note that $\theta_{\alpha,\alpha} = \theta_{\alpha}$.

Let $T(X,Y)=\nabla_XY-\nabla_YX-[X,Y]$ be the torsion tensor of
$\nabla$. We denote by the same letter the torsion tensor of type
(0,3) given by $T(X,Y,Z)=g(T(X,Y),Z)$. The Nijenhuis tensor is
expressed in terms of $\nabla$ as follows
\begin{eqnarray}\label{2}
N_{\alpha}(X,Y)  =  - \epsilon_{\alpha}4T_{\alpha}^{0,2}(X,Y)
                 + (\nabla_{J_{\alpha}X}J_{\alpha})(Y) -
(\nabla_{J_{\alpha}Y}J_{\alpha})(X) - (\nabla_Y
J_{\alpha})(J_{\alpha}X) +
(\nabla_XJ_{\alpha})(J_{\alpha}Y),\nonumber
\end{eqnarray}
where the (0,2)-part $T_{\alpha}^{0,2}$ of the torsion with
respect to $J_{\alpha}$ is given by
\begin{equation}\label{tr1}
T_{\alpha}^{0,2}(X,Y) = \frac{1}{4}\left( T(X,Y) +
\epsilon_{\alpha} T(J_{\alpha}X,J_{\alpha}Y) - \epsilon_{\alpha}
J_{\alpha}T(J_{\alpha}X,Y) - \epsilon_{\alpha}
J_{\alpha}T(X,J_{\alpha}Y)\right).
\end{equation}
A 3-form $\psi$ is of type (1,2)+(2,1) with respect to an almost
(para)complex structure $J_{\alpha}$ if and only if it satisfies
the equality \cite{I5}
\begin{equation}\label{3}
-
\epsilon_{\alpha}\psi(X,Y,Z)=\psi(J_{\alpha}X,J_{\alpha}Y,Z)+\psi(J_{\alpha}X,Y,J_{\alpha}Z)
+ \psi(X,J_{\alpha}Y,J_{\alpha}Z).
\end{equation}
{\bf Definition}. An almost paraquaternionic hermitian manifold
$(M,g,(H_{\alpha})\in \mathbb P)$ is a {\it PQKT manifold} if it
admits a metric paraquaternionic connection $\nabla$ with totally
skew symmetric torsion which is (1,2)+(2,1)-form with respect to
each $J_{\alpha}, \alpha=1,2,3$.  If the torsion 3-form is closed
then the manifold is said to be a {\it strong PQKT manifold}.

It follows that the holonomy group of $\nabla$ is a subgroup of
$Sp(n,\mathbb R)$.$Sp(1,\mathbb R)$.

By means of (\ref{1}), (\ref{2}) and (\ref{3}), the Nijenhuis
tensor $N_{\alpha}$ of $J_{\alpha},\alpha=1,2,3$, on a PQKT
manifold is given by
\begin{equation}\label{6}
N_{\alpha}(X,Y)  = - A_{\alpha}(Y)J_{\beta}X +
A_{\alpha}(X)J_{\beta}Y - J_{\alpha} A_{\alpha}(Y)J_{\gamma}X +
J_{\alpha} A_{\alpha}(X)J_{\gamma}Y,
\end{equation}
where
\begin{equation}\label{c1}
A_{\alpha}=\omega_{\beta} - \epsilon_{\alpha}J_{\alpha}
\omega_{\gamma}.
\end{equation}
\noindent {\bf Remark 1.} The torsion of $\nabla$ is a
(1,2)+(2,1)- form with respect to any (local) almost (para)complex
structure $J\in \mathbb P$ . In fact, it is sufficient that the
torsion is a (1,2)+(2,1)-form with respect to the only two almost
(para)complex structures of $(H)$ since Proposition 6.1 in
\cite{I5} gives the necessary expression of $N_{J_{\alpha}}$ by
$N_{J_{\beta}}$ and $N_{J_{\gamma}}$. Indeed, it is easy to see
that the formula from the Proposition 6.1 in \cite{I5} holds for
the (0,2)-part $T_{\alpha}^{0,2}, \alpha =1,2,3$, of the torsion.
Hence, the vanishing of the (0,2)-part of the torsion with respect
to any two almost (para)complex structures in $(H)$ implies the
vanishing of the (0,2)-part of $T$ with respect to the third one.

On a PQKT manifold there are three naturally associated 1-forms
defined by
\begin{equation}\label{n1}
t_{\alpha}(X)=\epsilon_{\alpha}
\frac{1}{2}\sum_{i=1}^{4n}\epsilon_{i} T(X,e_i,J_{\alpha}e_i),
\quad \alpha=1,2,3.
\end{equation}
Following \cite{I1}, we have
\begin{pro}\label{l1}
On a PQKT manifold  $J_1t_1=J_2t_2=J_3t_3.$
\end{pro}
{\it Proof.} Applying (\ref{3}) with respect to $J_{\beta}$ we
obtain
\begin{eqnarray}
t_{\alpha}(X) = \epsilon_{\alpha} \frac{1}{2}\sum_{i=1}^{4n}
\epsilon_{i} T(X,e_i,J_{\alpha}e_i) = - \frac{1}{2}\sum_{i=1}^{4n}
T(X,J_{\beta}e_i,J_{\gamma}e_i) \nonumber\\
               = \frac{1}{2}\sum_{i=1}^{4n}
               \epsilon_{i}T(J_{\beta}X,e_i,J_{\gamma}e_i)+\frac{1}{2}\sum_{i=1}^{4n}
               \epsilon_{i}T(J_{\beta}X,e_i,J_{\gamma}e_i)
            - \epsilon_{\alpha} \frac{1}{2} \sum_{i=1}^{4n} \epsilon_{i}
T(X,e_i,J_{\alpha}e_i).\nonumber
\end{eqnarray}
The last equality implies $t_{\alpha}=\epsilon_{\alpha}
\epsilon_{\beta}J_{\beta}t_{\gamma}$ which proves the assertion.
\hfill {\bf Q.E.D.}
\\
We introduce the torsion 1-form on PQKT manifolds by the equality
\begin{equation}\label{n1.1}
t(X)=-\frac{1}{2} \sum_{i=1}^{4n}
 \epsilon_{\alpha}\epsilon_{i}T(J_{\alpha}X,e_i,J_{\alpha}e_i).
\end{equation}
We need the following
\begin{lem}\label{l.10}
For a three form $T$ of type $(1,2)+(2,1)$ with respect to each
$J_{\alpha}$ one has
$$\sum_{i,j=1}^{4n}\epsilon_{i}\epsilon_{j}g(T(e_{i},e_{j}),T(J_{\gamma}e_{i},J_{\beta}e_{j}))=0,\sum_{i,j=1}^{4n}\epsilon_{i}\epsilon_{j}g(T(e_{i},e_{j}),T(J_{\beta}e_{i},J_{\beta}e_{j}))=-\frac{1}{3}\epsilon_{\beta}|T|^2,$$
where $|\cdot|^2$ denotes the norm with respect to the metric $g$.
\end{lem}
{\it Proof.} This proof is very similar to the proof of Lemma3.2
in \cite{I2} and we omit it. \hfill {\bf Q.E.D.}
\\
Further, we have
\begin{thm}\label{t0}
Let $(M,g,(J_{\alpha}\in \mathbb P)$ be a 4n-dimensional $PQKT$
manifold. Then the following identities hold
\begin{equation}\label{l1.1}
\sum_{i,j=1}^{4n}\epsilon_{i}(\nabla_{X}T)(J_{\alpha}Y,e_{i},J_{\alpha}e_{i})=-2\epsilon_{\alpha}(\nabla_{X}t)(Y);
\end{equation}
\begin{equation}\label{l1.2}
\sum_{i,j=1}^{4n}\epsilon_{i}\epsilon_{j}dT(e_{j},J_{\alpha}e_{j},e_{i}J_{\alpha}e_{i})=8\epsilon_{\alpha}\delta
t - 8\epsilon_{\alpha}|t|^2 + \frac{4}{3}\epsilon_{\alpha}|T|^2
,\sum_{i,j=1}^{4n}\epsilon_{i}\epsilon_{j}dT(e_{j},J_{\beta}e_{j},e_{i}J_{\gamma}e_{i})=0,
\end{equation}
where $\delta$ is the codifferential with respect to $g$.
\end{thm}
{\it Proof.} The formula $(\ref{l1.1})$ follows from $(\ref{1})$
and definition $(\ref{n1.1})$ of the torsion $1-$form by
straightforward calculations. To prove $(\ref{l1.2})$ we need the
expression of $dT$ in terms of $\nabla$ \cite{F2,I1},
\begin{equation}\label{l1.3}
dT(X,Y,Z,U)= {\sigma \atop XYZ}\left\{(\nabla_XT)(Y,Z,U) +
2g(T(X,Y),T(Z,U)\right\} - (\nabla_UT)(X,Y,Z),
\end{equation}
where $\sigma \atop XYZ$ denotes the cyclic sum of $X,Y,Z$. Taking
the appropriate trace in $(\ref{l1.3})$ and applying
Lemma$\ref{l.10}$ we obtain the first equality in $(\ref{l1.2})$.
Finally, from $(\ref{l1.3})$ combined with $(\ref{l1.1})$ and
Lemma$\ref{l.10}$ we get that
$$\sum_{i,j=1}^{4n}\epsilon_{i}\epsilon_{j}dT(e_{j},J_{\beta}e_{j},e_{i}J_{\gamma}e_{i})=-4\sum_{i,j=1}^{4n}\epsilon_{i}\epsilon_{j}g(T(e_{i},e_{j}),T(J_{\gamma}e_{i},J_{\beta}e_{j}))=0.$$
\hfill {\bf Q.E.D.}
\\
\begin{thm}\label{thm0}
Every PQKT is a paraquaternionic manifold.
\end{thm}
{\it Proof.} This is an immediate consequence of $(\ref{6})$ and
Theorem$\ref{thm0.4}$ \hfill {\bf Q.E.D.}

However, the converse to the above property is not always true. In
fact, we have

\begin{thm}\label{thm1}
Let $(M,g,(J_{\alpha}\in \mathbb P)$ be a 4n-dimensional ($n > 1$)
paraquaternionic manifold with $\mathbb P$-hermitian metric $g$.
Then $M$ admits a PQKT structure if and only if the following
conditions hold
\begin{equation}\label{4}
(d_{\alpha}F_{\alpha})^+ - (d_{\beta}F_{\beta})^+
=\frac{1}{2}\left(\epsilon_{\gamma}K_{\alpha}\wedge F_{\beta}-
\epsilon_{\beta}J_{\beta}K_{\beta}\wedge F_{\alpha}-
\epsilon_{\alpha}(K_{\beta}-J_{\alpha}K_{\alpha})\wedge
F_{\gamma}\right),
\end{equation}
where $(d_{\alpha}F_{\alpha})^+$ denotes the (1,2)+(2,1) part of
$(d_{\alpha}F_{\alpha})$ with respect to $J_{\alpha},
\alpha=1,2,3$ and the 1- forms $K_{\alpha}, \alpha=1,2,3$, are
given by
\begin{equation}\label{c2}
K_{\alpha}=
\frac{1}{1-n}\left(\epsilon_{\alpha}J_{\beta}\theta_{\alpha} +
\epsilon_{\beta}\theta_{\alpha,\gamma}\right).
\end{equation}
 The metric paraquaternionic
connection $\nabla$ with torsion 3-form of type (1,2)+(2,1) is
unique and is determined  by
\begin{equation}\label{5}
\nabla =\nabla^g +\frac{1}{2}\left((d_{\alpha}F_{\alpha})^+
-\frac{1}{2} \left(\epsilon_{\alpha}J_{\alpha}K_{\alpha}\wedge
F_{\gamma}+\epsilon_{\gamma}K_{\alpha}\wedge
F_{\beta}\right)\right),
\end{equation}
where $\nabla^g$ is the Levi-Civita connection of $g$.
\end{thm}
{\it Proof.} To prove the 'if' part, let $\nabla$ be a metric
paraquaternionic connection satisfying (\ref{1}) which torsion $T$
has the required properties. We follow the scheme in \cite{I1}.
Since $T$ is skew-symmetric, we have
\begin{equation}\label{5'}
\nabla =\nabla^g + \frac{1}{2}T.
\end{equation}
We obtain using (\ref{1}) and (\ref{5'}) that
\begin{eqnarray}\label{7}
\frac{1}{2}\left(T(X,J_{\alpha}Y,Z) + (T(X,Y,J_{\alpha}Z)\right) =
-g\left((\nabla_X^gJ_{\alpha})Y,Z\right)\\ -
\omega_{\beta}(X)F_{\gamma}(Y,Z) - \epsilon_{\gamma}
\omega_{\gamma}(X)F_{\beta}(Y,Z).\nonumber
\end{eqnarray}
The tensor $\nabla^gJ_{\alpha}$ is decomposed by parts according
to
 $\nabla J_{\alpha} =
(\nabla J_{\alpha})^{2,0} + (\nabla J_{\alpha})^{0,2}$, where
\cite{I5,G2}
\begin{equation}\label{10}
g\left((\nabla_X^gJ_{\alpha})^{2,0}Y,Z\right) = -
\frac{1}{2}\left(\epsilon_{\alpha}(dF_{\alpha})^+(X,J_{\alpha}Y,J_{\alpha}Z)
+ (dF_{\alpha})^+(X,Y,Z)\right)
\end{equation}
\begin{equation}\label{9}
g\left((\nabla_X^gJ_{\alpha})^{0,2}Y,Z\right) =
\frac{1}{2}\left(g(N_{\alpha}(X,Y),J_{\alpha}Z) -
g(N_{\alpha}(X,Z),J_{\alpha}Y)-
g(N_{\alpha}(Y,Z),J_{\alpha}X)\right)
\end{equation}
Taking the (2,0) part in (\ref{7}) we obtain using (\ref{10}) that
\begin{eqnarray}\label{i2}
T(X,J_{\alpha}Y,Z) + T(X,Y,J_{\alpha}Z) =
(\epsilon_{\alpha}dF_{\alpha}^+(X,J_{\alpha}Y,J_{\alpha}Z)+
dF_{\alpha}^+(X,Y,Z))\\ -
C_{\alpha}(X)F_{\gamma}(Y,Z)+\epsilon_{\gamma}C_{\alpha}(J_{\alpha}X)F_{\beta}(Y,Z),\nonumber
\end{eqnarray}
where
\begin{equation}\label{c3}
C_{\alpha}=\omega_{\beta}+\epsilon_{\alpha}J_{\alpha}\omega_{\gamma}.
\end{equation}
The cyclic sum of (\ref{i2}) and  the fact that $T$ and
$(dF_{\alpha})^+$ are (1,2)+(2,1)-forms with respect to each
$J_{\alpha}$, gives
\begin{equation}\label{i3}
T= (d_{\alpha}F_{\alpha})^+ -\frac{1}{2}
\left(\epsilon_{\alpha}J_{\alpha}C_{\alpha}\wedge
F_{\gamma}+\epsilon_{\gamma}C_{\alpha}\wedge F_{\beta}\right).
\end{equation}
Further, we take the contractions in (\ref{i3}) to get
\begin{eqnarray}\label{c4}
J_{\alpha}t_{\alpha} = - \theta_{\alpha} -
\epsilon_{\gamma}J_{\beta}C_{\alpha}, \nonumber\\
J_{\alpha}t_{\alpha} =
\epsilon_{\beta}J_{\gamma}\theta_{\beta,\alpha} - n
\epsilon_{\alpha}J_{\gamma}C_{\beta},\\ J_{\alpha}t_{\alpha} = -
\epsilon_{\gamma}J_{\beta} \theta_{\gamma,\alpha} - n
\epsilon_{\beta}J_{\alpha}C_{\gamma}\nonumber
\end{eqnarray}
Using Proposition~\ref{l1}, (\ref{c1}) and (\ref{c3}), we obtain
consequently from (\ref{c4})  that
\begin{equation}\label{c5}
\epsilon_{\alpha}A_{\alpha}= -J_{\alpha}C_{\beta}
+\epsilon_{\alpha}J_{\gamma}C_{\gamma} =
J_{\beta}\left(\theta_{\gamma} - \theta_{\beta}\right),
\end{equation}
\begin{equation}\label{c6}
(n-1)\epsilon_{\gamma}J_{\beta}C_{\alpha}=  \theta_{\alpha} +
\epsilon_{\alpha}J_{\beta}\theta_{\alpha,\gamma}.
\end{equation}
Then (\ref{4}) and (\ref{c2}) follow from (\ref{i3}) and
(\ref{c6}).

For the converse, we define $\nabla$ by (\ref{5}). To complete the
proof we have to show that $\nabla$ is a paraquaternionic
connection. We calculate
\begin{eqnarray}
g\left((\nabla_XJ_{\alpha})Y,Z\right) =
g\left((\nabla_X^gJ_{\alpha})Y,Z\right)+
\frac{1}{2}\left(T(X,J_{\alpha}Y,Z) + T(X,Y,J_{\alpha}Z)\right)
\nonumber\\
                                      =
- \omega_{\beta}(X)F_{\gamma}(Y,Z) -
\epsilon_{\gamma}\omega_{\gamma}(X)F_{\beta}(Y,Z),\nonumber
\end{eqnarray}
where we used (\ref{10}), (\ref{9}), (\ref{c5}), (\ref{c2}),
(\ref{c1}), (\ref{c3}) and the compatibility condition (\ref{4})
to get the last equality. The uniqueness of $\nabla $ follows from
(\ref{5}). \hfill {\bf Q.E.D.}

In the case of hyper-paraK\"ahler manifold with torsion (briefly
HPKT), $K_{\alpha}=dF_{\alpha}^-=0$ and Theorem~\ref{thm1} is a
consequence of the general results in \cite{I5} which imply that
on a para-Hermitian manifold there exists a unique linear
connection with totally skew-symmetric torsion preserving the
metric and the (para)complex structure.
\\
As a consequence of the proof of Theorem~\ref{thm1}, we get
\begin{pro}\label{cc1}
The Nijenhuis tensors  of a PQKT manifold depend only on the
difference between the Lie forms. In particular, the almost
(para)complex structures $J_{\alpha}$ on a PQKT manifold
$(M,(J_{\alpha})\in \mathbb P,g,\nabla )$ are integrable if and
only if
$$\theta_{\alpha}=\theta_{\beta}= \theta_{\gamma}$$
\end{pro}
{\it Proof.} The Nijenhuis tensors  are given by (\ref{6}) and
(\ref{c5}). \hfill {\bf Q.E.D.}

\begin{co}\label{tt1}
On a 4n-dimensional PQKT manifold the following formulas hold
$$J_{\beta}\theta_{\alpha,\gamma}= -
J_{\gamma}\theta_{\alpha,\beta},$$
\begin{equation}\label{per1}
(n^2+n)\theta_{\alpha} -n\theta_{\beta} - n^2\theta_{\gamma} -
\epsilon_{\beta}J_{\gamma}\theta_{\beta,\alpha} -
n\epsilon_{\gamma}J_{\alpha}\theta_{\gamma,\beta} +
(n+1)\epsilon_{\alpha}J_{\beta}\theta_{\alpha,\gamma}=0.
\end{equation}
If $n=1$ then  \hspace{4cm} $\theta_{\alpha} = - \epsilon_{\alpha}
J_{\beta} \theta_{\alpha,\gamma} =
\epsilon_{\alpha}J_{\gamma}\theta_{\alpha,\beta}$.
\end{co}
{\it Proof.} The first formula follows directly from the system
(\ref{c4}). Solving the system (\ref{c4}) with respect to
$C_{\alpha}$ we obtain
\begin{equation}\label{per2}
(n^3-1)\epsilon_{\gamma}J_{\beta}C_{\alpha} =
(\theta_{\alpha}+\epsilon_{\beta}J_{\gamma}\theta_{\beta,\alpha})
+
n(\theta_{\beta}+\epsilon_{\gamma}J_{\alpha}\theta_{\gamma,\beta})
+
n^2(\theta_{\gamma}+\epsilon_{\alpha}J_{\beta}\theta_{\alpha,\gamma}).
\end{equation}
Then (\ref{per1}) is a consequence of (\ref{per2}) and (\ref{c6}).
 The last assertion follows from (\ref{c6}) .
\hfill {\bf Q.E.D.}
\begin{co}\label{cc2}
On a 4n-dimensional ($n>1$) PQKT manifold the $sp(1,\mathbb
R)$-connection 1-forms are given by
\begin{equation}\label{c7}
\omega_{\beta} =
\frac{1}{2}\epsilon_{\alpha}J_{\beta}\left(\theta_{\gamma} -
\theta_{\beta} + \frac{1}{1-n}\theta_{\alpha}\right) +
\frac{1}{2(1-n)}\epsilon_{\beta}\theta_{\alpha,\gamma}.
\end{equation}
\end{co}
{\it Proof.} The proof follows in a straightforward way from
(\ref{c5}), (\ref{c6}), (\ref{c1}) and (\ref{c3}). \hfill {\bf
Q.E.D.}

Theorem~\ref{thm1} and the above formulas lead to the following
criterion
\begin{pro}\label{cc3}
Let $(M,g,(J_{\alpha}))$ be a 4n-dimensional ($n>1$) PQKT
manifold. The following conditions are equivalent:

i) $(M,g,(H))$ is a local HPKT manifold;

ii) $d_{\alpha}F_{\alpha}^+=d_{\beta}F_{\beta}^+=
d_{\gamma}F_{\gamma}^+$;

iii) $ \theta_{\alpha} =
-\epsilon_{\alpha}J_{\beta}\theta_{\alpha,\gamma}. $
\end{pro}
{\it Proof.} If $(M,g,(H))$ is a HPKT manifold, the connection
1-forms $\omega_{\alpha}=0, \alpha=1,2,3$. Then ii) and iii)
follow from (\ref{c3}), (\ref{c6}), (\ref{c2}) and (\ref{4}).

If iii) holds, then (\ref{c6}) and (\ref{c5}) yield
$C_{\alpha}=A_{\alpha}=0, \alpha=1,2,3$, since $n>1$.
Consequently, $2\omega_{\alpha}= J_{\beta}C_{\beta}-
J_{\beta}A_{\beta}=0$ by (\ref{c3}) and (\ref{c1}). Thus the
equivalence of i) and iii) is proved.

Let ii) holds. Then we compute that $\theta_{\alpha}=
J_{\gamma}\theta_{\beta,\alpha}$. Since $n>1$, the equality
(\ref{per2}) leads to $C_{\alpha}=0, \alpha=1,2,3$, which forces
$\omega_{\alpha}=0, \alpha=1,2,3$  as above. This completes the
proof. \hfill {\bf Q.E.D.}

The next theorem shows that PQKT manifolds are stable under a
conformal transformations.
\begin{thm}\label{thm3}
Let $(M,g,(J_{\alpha}),\nabla)$ be a 4n-dimensional PQKT manifold.
Then every pseudo-Riemannian metric $\bar{g}$ in the conformal
class [g] admits a PQKT connection. If $\bar{g}=fg$ for a function
$f$ then the PQKT connection $\bar{\nabla}$ corresponding to
$\bar{g}$
 is given by
\begin{eqnarray}\label{z1}
\bar{g}(\bar{\nabla}_XY,Z)=fg(\nabla_XY,Z) +
\frac{1}{2}\left(df(X)g(Y,Z) + df(Y)g(X,Z) - df(Z)g(X,Y)\right)\\
- \frac{1}{2}\left(\epsilon_{\alpha}J_{\alpha}df\wedge F_{\alpha}
+\epsilon_{\beta}J_{\beta}df\wedge F_{\beta}+
\epsilon_{\gamma}J_{\gamma}df\wedge
F_{\gamma}\right)(X,Y,Z).\nonumber
\end{eqnarray}
The torsion tensors $T$ and $\bar{T}$ and the torsion 1-forms $t$
and $\bar{t}$ of $\nabla$ and $\bar{\nabla}$ are related by
\begin{equation}\label{z4}
\bar{T}= fT - \epsilon_{\alpha}J_{\alpha}df\wedge F_{\alpha} -
\epsilon_{\beta}J_{\beta}df\wedge F_{\beta}-
\epsilon_{\gamma}J_{\gamma}df\wedge F_{\gamma}.
\end{equation}
\begin{equation}\label{z5}
\bar{t}= t - (2n+1)d\ln f.
\end{equation}
\end{thm}

{\it Proof.} First we assume $n>1$. We apply Theorem~\ref{thm1} to
the paraquaternionic Hermitian manifold
$(M,\bar{g}=fg,(J_{\alpha})\in \mathbb P)$. We denote the objects
corresponding to the metric $\bar{g}$ by a line above the symbol
e.g. $\bar{F_{\alpha}}$ denotes the K\"ahler form of $J_{\alpha}$
with respect to $\bar{g}$. An easy calculation gives the following
sequence of formulas
\begin{equation}\label{z2}
d_{\alpha}\bar{F}_{\alpha}^+ = -
\epsilon_{\alpha}J_{\alpha}df\wedge F_{\alpha} +
fd_{\alpha}F_{\alpha}^+; \quad \bar{\theta}_{\alpha} =
\theta_{\alpha} +(2n-1)d\ln f; \quad \bar{\theta}_{\alpha,\gamma}
= \theta_{\alpha,\gamma} + \epsilon_{\gamma}J_{\beta}d\ln f.
\end{equation}
We substitute (\ref{z2}) into (\ref{c2}), (\ref{c5}) and
(\ref{c7}) to get
\begin{equation}\label{z3}
\bar{K}_{\alpha} = K_{\alpha} -2\epsilon_{\alpha}J_{\beta}d\ln f,
\quad \bar{A}=A, \quad \bar{\omega}_{\alpha} = \omega_{\alpha}
-\epsilon_{\alpha}J_{\beta}d\ln f.
\end{equation}
Using (\ref{z2}) and (\ref{z3}) we verify that the conditions
(\ref{4}) with respect to the metric $\bar{g}$ are fulfilled.
Theorem~\ref{thm1} implies that there exists a PQKT connection
$\bar{\nabla}$ with respect to  $(\bar{g},P)$. Using the well
known relation between the Levi-Civita connections of conformally
equivalent metrics, (\ref{z2}) and (\ref{z3}), we obtain
(\ref{z1}) from (\ref{5}).
\\
Using (\ref{z1}), we get (\ref{z4}) and consequently (\ref{z5}).
\hfill {\bf Q.E.D.}
\\
Namely, any conformal metric of a PQK, HPK or HPKT manifold will
give a PQKT manifold. This leads to the notion of {\it locally
conformally PQK (resp. locally conformally HPK, resp. locally
conformally HPKT) manifolds} (briefly l.c.PQK (resp. l.c.HPK,
resp. l.c.HPKT) manifolds) in the context of PQKT geometry.
\\
We recall that a paraquaternionic Hermitian manifold $(M,g,\mathbb
P)$ is said to be  l.c.PQK (resp. l.c.HPK, resp. l.c.HPKT)
manifold if each point $p\in M$ has a neighborhood $U_p$ such that
$g \Big|_{U_p}$ is conformally equivalent to a PQK (resp.HPK,
resp.HPKT) metric.
\\
For example, the Kodaira-Thurston surface modeled on $\widetilde
{SL(2,\mathbb R)}$$\times \mathbb R/$$\Gamma$ is on example of a
compact l.c.HPKT which is not globaly conformal HPKT \cite{I5}.
\\
Theorem~\ref{thm3},  Theorem~\ref{thm1} together with
Proposition~\ref{cc1} and Proposition~\ref{cc3} imply  the
following
\begin{co}\label{u1}
Every l.c.PQK manifold admits a PQKT structure.

Further, if  $(M,g,(J_{\alpha}),\nabla)$ is a 4n-dimensional $n>1$
PQKT manifold then:

i) $(M,g,(J_{\alpha}),\nabla)$ is a l.c.PQK manifold if and only
if
\begin{equation}\label{u3}
T = - \frac{1}{2n+1}\left(t_{\alpha}\wedge F_{\alpha} +
t_{\beta}\wedge F_{\beta}+t_{\gamma}\wedge F_{\gamma}\right),
\quad dt=0;
\end{equation}
\indent ii) $(M,g,(J_{\alpha}),\nabla)$ is a l.c.HPKT manifold if
and only if the 1-form
$\theta_{\alpha}+\epsilon_{\alpha}J_{\beta}\theta_{\alpha,\gamma}$
is closed i.e.
$$ d(\theta_{\alpha}+\epsilon_{\alpha}J_{\beta}\theta_{\alpha,\gamma})=0; $$
\indent iii)  $(M,g,(J_{\alpha}),\nabla)$ is a l.c.HPK manifold if
an only if (\ref{u3}) holds and $$
\theta_{\alpha}+\epsilon_{\alpha}J_{\beta}\theta_{\alpha,\gamma}=\frac{2(1-n)}{2n+1}t.
$$
\end{co}

\section{Curvature of a PQKT space}

Let $R=[\nabla,\nabla]-\nabla_{[,]}$ be the curvature tensor of
type (1,3) of $\nabla$. We denote the curvature tensor of type
(0,4) $R(X,Y,Z,V)=g(R(X,Y)Z,V)$ by the same letter. There are
three Ricci forms and three scalar function given by $$
\rho_{\alpha}(X,Y)=\frac{1}{2}
\sum_{i=1}^{4n}\epsilon_{i}R(X,Y,e_i,J_{\alpha}e_i), \quad
\alpha=1,2,3, $$ $$ Scal_{\alpha,\beta}=
\sum_{i=1}^{4n}\epsilon_{i}
\epsilon_{\alpha}\rho_{\alpha}(e_i,J_{\beta}e_i), \quad
\alpha,\beta=1,2,3. $$
\begin{pro}\label{p1}
The curvature of a PQKT manifold $(M,g,(J_{\alpha}),\nabla)$
satisfies the following relations
\begin{equation}\label{11}
[R(X,Y),J_{\alpha}] = \frac{1}{n}\left(-
\epsilon_{\alpha}\rho_{\gamma}(X,Y) \otimes J_{\beta} +
\epsilon_{\alpha}\rho_{\beta}(X,Y) \otimes J_{\gamma}\right),
\end{equation}
\begin{equation}\label{12}
\epsilon_{\gamma} \rho_{\alpha} =d\omega_{\alpha} +
\epsilon_{\alpha} \omega_{\beta}\wedge\omega_{\gamma}.
\end{equation}
\end{pro}
{\it Proof.} We follow the classical scheme (see e.g. \cite{I5}).
Using (\ref{1}), we obtain $$ [R(X,Y),J_{\alpha}]=
(d\omega_{\beta} +
\epsilon_{\beta}\omega_{\gamma}\wedge\omega_{\alpha})(X,Y) \otimes
J_{\gamma} + \epsilon_{\gamma} (d\omega_{\gamma} +
\epsilon_{\gamma}\omega_{\alpha}\wedge\omega_{\beta})(X,Y) \otimes
J_{\beta}.
$$ Taking the trace in the last equality, we get
\begin{eqnarray}
\rho_{\alpha}(X,Y) = \frac{1}{2}
\sum_{i=1}^{4n}\epsilon_{i}R(X,Y,e_i,J_{\alpha}e_i)=- \frac{1}{2}
\sum_{i=1}^{4n}\epsilon_{\beta}\epsilon_{i}R(X,Y,J_{\beta}e_i,J_{\alpha}J_{\beta}e_i)
\nonumber\\
                   =
-\frac{1}{2} \sum_{i=1}^{4n}\epsilon_{i}R(X,Y,e_i,J_{\alpha}e_i)
+2n\epsilon_{\gamma}(d\omega_{\alpha} + \epsilon_{\alpha}
\omega_{\beta}\wedge\omega_{\gamma})(X,Y).\nonumber
\end{eqnarray}
\hfill {\bf Q.E.D.}

Using Proposition~\ref{p1} we find a simple necessary and
sufficient condition a PQKT manifold to be a HPKT one, i.e. the
holonomy group of $\nabla$ to be a subgroup of $Sp(n,\mathbb R)$.
\begin{pro}\label{p2}
A 4n-dimensional $(n>1)$ PQKT manifold is a local HPKT manifold if
and only if all the three Ricci forms vanish, i.e
$\rho_1=\rho_2=\rho_3=0$.
\end{pro}
{\it Proof.} If a PQKT manifold is a HPKT manifold then the
holonomy group of $\nabla$ is contained in $Sp(n,\mathbb R)$. This
implies $\rho_{\alpha}=0, \quad \alpha=1,2,3$.

For the converse, let the three Ricci forms vanish. The equations
(\ref{12}) mean that the curvature of the $Sp(1,\mathbb R)$
connection on $\mathbb P$ vanish. Then there exists a local basis
$(I_{\alpha}, \alpha =1,2,3)$ of almost (para)complex structures
on $\mathbb P$ and each $I_{\alpha}$ is $\nabla$-parallel i.e. the
corresponding connection 1-forms $\omega_{I_{\alpha}}=0, \alpha
=1,2,3$. Then each $I_{\alpha}$ is a (para)complex structure, by
(\ref{6}) and (\ref{c1}). This implies that the PQKT manifold is a
local HPKT manifold. \hfill {\bf Q.E.D.}

The Ricci tensor $Ric$ and scalar curvatures $Scal$ and
$Scal_{\alpha}$ of the PQKT connection $\nabla$ are defined by
$$Ric(X,Y)=\sum_{i=1}^{4n}\epsilon_{i}R(e_i,X,Y,e_i), Scal=\sum_{i=1}^{4n}\epsilon_{i}Ric(e_i,e_i), Scal_{\alpha}=-\sum_{i=1}^{4n}\epsilon_{i}Ric(e_i,J_{\alpha}e_i).$$
We denote by $Ric^g,$$Scal^g,$$\rho_{\alpha}^g,$ etc. the
corresponding objects for the metric $g,$ i.e. the same objects
taken with respect to the Levi-Civita connection $\nabla^g.$We may
consider $(g,J_{\alpha})$ as an almost (para)Hermitian structure.
The tensor
$\rho_{\alpha}^{\star}(X,Y)=\rho_{\alpha}^{g}(X,J_{\alpha}Y)$ is
known as the $\star -$ Ricci tensor of the almost (para)Hermitian
structure. It is equal to $\rho_{\alpha}^{\star}(X,Y)=-
\sum_{i=1}^{2n}R^g(e_i,X,J_{\alpha}Y,J_{\alpha}e_i)$ by the
Bianchi identity. The function $Scal_{\alpha}^g$ is known also as
the $\star-$scalar curvature. In general, the $\star -$Ricci
tensor is not symmetric and the $\star -$ Einstein condition is a
strong condition. We shall see in this section that the scalar
curvature functions are not independent and we define a new scalar
invariant, the "paraquaternionic $\star -$ scalar curvature" of a
PQKT space.
\\
Our main technical result is the following
\begin{pro}\label{tir}
Let $(M,g,(J_{\alpha}),\nabla)$ be a 4n-dimensional PQKT manifold.
The following formulas hold
\begin{eqnarray}\label{ti20}
 n \epsilon_{\alpha} \rho_{\alpha}(X,J_{\alpha}Y)+
\epsilon_{\beta} \rho_{\beta}(X,J_{\beta}Y) + \epsilon_{\gamma}
\rho_{\gamma}(X,J_{\gamma}Y)= \\                      nRic(X,Y)
+\frac{n}{4}\epsilon_{\alpha} (dT)_{\alpha}(X,J_{\alpha}Y)-
n(\nabla_X t)(Y);\nonumber
\end{eqnarray}
\begin{equation}\label{ti22}
 (n-1) \epsilon_{\alpha} \rho_{\alpha}(X,J_{\alpha}Y)=
\frac{n(n-1)}{n+2}Ric(X,Y) - \frac{n(n-1)}{(n+2)}(\nabla_Xt)Y+
\end{equation}
$$+ \frac{n}{4(n+2)}\left\{(n+1)\epsilon_{\alpha}
(dT)_{\alpha}(X,J_{\alpha}Y)- \epsilon_{\beta}
(dT)_{\beta}(X,J_{\beta}Y) -\epsilon_{\gamma}
(dT)_{\gamma}(X,J_{\gamma}Y)\right\},$$
where
$(dT)_{\alpha}(X,Y)=\sum_{i=1}^{4n}\epsilon_{i}dT(X,Y,e_i,J_{\alpha}e_i)$.
\end{pro}
{\it Proof.} Since the torsion is a 3-form, we have \cite{F2,I1}
\begin{equation}\label{sof}
(\nabla^g_XT)(Y,Z,U) = (\nabla_XT)(Y,Z,U) + \frac{1}{2} {\sigma
\atop XYZ} \left\{g(T(X,Y),T(Z,U)\right\},
\end{equation}
where ${\sigma \atop XYZ}$ denote the cyclic sum of $X,Y,Z$.

The exterior derivative $dT$ is given by
\begin{eqnarray}\label{13}
dT(X,Y,Z,U)= {\sigma \atop XYZ}\left\{(\nabla_XT)(Y,Z,U) +
g(T(X,Y),T(Z,U)\right\}\\
           - (\nabla_UT)(X,Y,Z) + {\sigma \atop XYZ}
\left\{g(T(X,Y),T(Z,U)\right\}. \nonumber
\end{eqnarray}
The first Bianchi identity for $\nabla$ states
\begin{equation}\label{14}
{\sigma \atop XYZ}R(X,Y,Z,U)= {\sigma \atop
XYZ}\left\{(\nabla_XT)(Y,Z,U) + g(T(X,Y),T(Z,U)\right\}.
\end{equation}
We denote by $B$ the Bianchi projector i.e. $B(X,Y,Z,U)={\sigma
\atop XYZ}R(X,Y,Z,U)$.

The curvature $R^g$ of the Levi-Civita connection is connected by
$R$ in the following way
\begin{eqnarray}\label{15}
R^g(X,Y,Z,U) = R(X,Y,Z,U) - \frac{1}{2} (\nabla_XT)(Y,Z,U)
+\frac{1}{2} (\nabla_YT)(X,Z,U)\nonumber \\
             -
\frac{1}{2}g(T(X,Y),T(Z,U)) - \frac{1}{4}g(T(Y,Z),T(X,U)) -
\frac{1}{4}g(T(Z,X),T(Y,U)).
\end{eqnarray}
Define $D$ by $D(X,Y,Z,U)=R(X,Y,Z,U)-R(Z,U,X,Y)$, we obtain from
(\ref{15})
\begin{eqnarray}\label{tir1}
 D(X,Y,Z,U) =\\
\frac{1}{2} (\nabla_XT)(Y,Z,U)-\frac{1}{2}(\nabla_YT)(X,Z,U) -
\frac{1}{2}(\nabla_ZT)(U,X,Y)+
\frac{1}{2}(\nabla_UT)(Z,X,Y),\nonumber
\end{eqnarray}
since $D^g$ of $R^g$ is zero.

Using (\ref{11}) and (\ref{14}), we find the following relation
between the Ricci tensor and the Ricci forms
\begin{eqnarray}\label{16}
\rho_{\alpha}(X,Y) = -\frac{1}{2} \sum_{i=1}^{4n}\left
(\epsilon_{i}(R(Y,e_i,X,J_{\alpha}e_i)+
\epsilon_{i}R(e_i,X,Y,J_{\alpha}e_i)\right) +\\
      + \frac{1}{2}
\sum_{i=1}^{4n}\epsilon_{i}B(X,Y,e_i,J_{\alpha}e_i)=
-\frac{1}{2}Ric(Y,J_{\alpha}X) +\frac{1}{2}Ric(X,J_{\alpha}Y)
+\frac{1}{2}\sum_{i=1}^{4n}\epsilon_{i}B(X,Y,e_i,J_{\alpha}e_i)
\nonumber\\
                   +
\frac{1}{2n}\left\{- \epsilon_{\alpha} \rho_{\beta}(J_{\gamma}Y,X)
+\epsilon_{\alpha}\rho_{\beta}(J_{\gamma}X,Y) - \epsilon_{\alpha}
\rho_{\gamma}(J_{\beta}X,Y)+ \epsilon_{\alpha}
\rho_{\gamma}(J_{\beta}Y,X)\right\}. \nonumber
\end{eqnarray}

On the other hand,  using  (\ref{11}), we calculate
\begin{eqnarray}\label{19}
\sum_{i=1}^4\epsilon_{i}D(X,e_i,J_{\alpha}e_i,Y)
=\sum_{i=1}^{4n}\{\epsilon_{i}R(X,e_i,J_{\alpha}e_i,Y)+
\epsilon_{i}R(Y,e_i,J_{\alpha}e_iX)\}\\
  =Ric(Y,J_{\alpha}X)+Ric(X,J_{\alpha}Y)\nonumber \\
  +
\frac{1}{n}\left\{\epsilon_{\alpha} \rho_{\beta}(X,J_{\gamma}Y)
+\epsilon_{\alpha} \rho_{\beta}(Y,J_{\gamma}X) -\epsilon_{\alpha}
\rho_{\gamma}(Y,J_{\beta}X)-\epsilon_{\alpha}
\rho_{\gamma}(X,J_{\beta}Y)\right\}. \nonumber
\end{eqnarray}
Combining (\ref{16}) and (\ref{19}), we derive
\begin{equation}\label{20}
n\epsilon_{\alpha} \rho_{\alpha}(X,J_{\alpha}Y)+ \epsilon_{\beta}
\rho_{\beta}(X,J_{\beta}Y) + \epsilon_{\gamma}
\rho_{\gamma}(X,J_{\gamma}Y)=
\end{equation}
$$
  nRic(X,Y)+\frac{n}{2}\epsilon_{\alpha}B_{\alpha}(X,J_{\alpha}Y)+
\frac{n}{2}\epsilon_{\alpha}D_{\alpha}(X,J_{\alpha}Y), $$ where
the tensors $B_{\alpha}$ and $D_{\alpha}$  are  defined by
$B_{\alpha}(X,Y)=\sum_{i=1}^{4n}\epsilon_{i}
B(X,Y,e_i,J_{\alpha}e_i)$ and \\
$ D_{\alpha}(X,Y)=\sum_{i=1}^{4n}\epsilon_{i}
D(X,e_i,J_{\alpha}e_i,Y) $. Taking into account (\ref{tir1}), we
get the expression
\begin{equation}\label{tir2}
D_{\alpha}(X,Y)=-(\nabla_Xt)(J_{\alpha}Y)
-(\nabla_Yt)(J_{\alpha}X)\quad \alpha=1,2,3.
\end{equation}
To calculate $B_{\alpha}+D_{\alpha}$ we use (\ref{13}) twice and
(\ref{tir2}). After some calculations, we derive
\begin{equation}\label{tir4}
B_{\alpha}(X,Y)+ D_{\alpha}(X,Y) =\frac{1}{2}  \sum_{i=1}^{4n}
\epsilon_{i}dT(X,Y,e_i,J_{\alpha}e_i) - 2(\nabla_Xt)(J_{\alpha}Y),
\quad
\\
\alpha=1,2,3.
\end{equation}
We substitute (\ref{tir4}) into (\ref{20}). Solving  the obtained
system, we obtain
\begin{eqnarray}\label{21}
(n-1)(\epsilon_{\alpha} \rho_{\alpha}(X,J_{\alpha}Y)-
\epsilon_{\beta} \rho_{\beta}(X,J_{\beta}Y)) =
\frac{n}{2}(\epsilon_{\alpha}(dT)_{\alpha}(X,J_{\alpha}Y) -
\epsilon_{\beta}(dT)_{\beta}(X,J_{\beta}Y)).
\end{eqnarray}
Finally, (\ref{20}) and (\ref{tir4}) imply (\ref{ti22}). \hfill
{\bf Q.E.D.}
\\
\begin{thm}\label{thm10}
On a PQKT manifold $(M^{4n},g,(J_{\alpha})\in \mathbb P)$ $(n>1)$
the (2,0)+(0,2)-parts of the Ricci forms $\rho_{\alpha}$,
$\rho_{\beta}$ with respect to $J_{\gamma}$ coincide in sense that
the next identity holds
\begin{equation}\label{21.2}
\rho_{\alpha}(J_{\beta}X,J_{\beta}Y)+\epsilon_{\beta}\rho_{\alpha}(X,Y)+\epsilon_{\gamma}\rho_{\gamma}(J_{\beta}X,Y)+\epsilon_{\gamma}\rho_{\gamma}(X,J_{\beta}Y)=0
\end{equation}
\end{thm}
{\it Proof}. We need the following
\begin{lem}\label{21.1}
The tensors $L_{\alpha}(X,Y)= \sum_{i=1}^{4n}
\epsilon_{i}g(T(X,e_i),T(Y,J_{\alpha}e_i),$ for $\alpha=1,2,3$ are
related by:
\begin{equation}\label{e21.1}
L_{\alpha}(J_{\beta}X,J_{\beta}Y)+
\epsilon_{\beta}L_{\alpha}(X,Y)+
\epsilon_{\gamma}L_{\gamma}(J_{\beta}X,Y)+
\epsilon_{\gamma}L_{\gamma}(X,J_{\beta}Y)=0
\end{equation}
\end{lem}
{\it Proof.} The formula $(\ref{e21.1})$ follows from equalities
$$L_{\alpha}(X,Y)= \sum_{i=1}^{4n}
\epsilon_{i}g(T(X,e_i),T(Y,J_{\alpha}e_i))= - \sum_{i=1}^{4n}
\epsilon_{i}g(T(X,J_{\alpha}e_i),T(Y,e_i))=$$
$$=\sum_{i,j=1}^{4n} \epsilon_{i}\epsilon_{j}T(X,e_i,e_j)T(e_j,Y,J_{\alpha}e_i)=\sum_{i,j=1}^{4n} \epsilon_{i}\epsilon_{j}T(X,e_i,e_j)g(T(e_j,Y),J_{\alpha}e_i)=$$
$$=\sum_{i,j=1}^{4n} \epsilon_{i}\epsilon_{j}T(X,e_j,e_i)g(e_i,J_{\alpha}T(Y,e_j))= - \sum_{j=1}^{4n} \epsilon_{j}g(T(X,e_j),J_{\alpha}T(Y,e_j))$$
and property $(\ref{3})$. \hfill {\bf Q.E.D.}
\\
From the first Bianchi identity the next sequence of equalities
\begin{gather*}
2\rho_{\alpha}(J_{\beta}X,J_{\beta}Y)+\sum_{i=1}^{4n}\epsilon_{i}R(J_{\beta}Y,e_i,J_{\beta}X,J_{\alpha}e_i)+\sum_{i=1}^{4n}\epsilon_{i}R(e_i,J_{\beta}X,J_{\beta}Y,J_{\alpha}e_i)=2\epsilon_{\gamma}(\nabla_{J_{\beta}X}t)J_{\gamma}Y-\\
2\epsilon_{\gamma}(\nabla_{J_{\beta}Y}t)J_{\gamma}X+\sum_{i=1}^{4n}\epsilon_{i}(\nabla_{e_{i}}T)(J_{\beta}X,J_{\beta}Y,J_{\alpha}e_i)+\sum_{i=1}^{4n}\epsilon_{i}g(T(J_{\beta}X,J_{\beta}Y),T(e_i,J_{\alpha}e_i))-2L_{\alpha}(J_{\beta}X,J_{\beta}Y).
\end{gather*}
\begin{gather*}
2\epsilon_{\beta}\rho_{\alpha}(X,Y)+\sum_{i=1}^{4n}\epsilon_{i}\epsilon_{\beta}R(Y,e_i,X,J_{\alpha}e_i)+\sum_{i=1}^{4n}\epsilon_{i}\epsilon_{\beta}R(e_i,X,Y,J_{\alpha}e_i)=-2\epsilon_{\beta}(\nabla_{X}t)J_{\alpha}Y+\\
2\epsilon_{\beta}(\nabla_{Y}t)J_{\alpha}X+\sum_{i=1}^{4n}\epsilon_{i}\epsilon_{\beta}(\nabla_{e_{i}}T)(X,Y,J_{\alpha}e_i)+\sum_{i=1}^{4n}\epsilon_{i}\epsilon_{\beta}g(T(X,Y),T(e_i,J_{\alpha}e_i))-2\epsilon_{\beta}L_{\alpha}(X,Y).
\end{gather*}
\begin{gather*}
2\epsilon_{\gamma}\rho_{\gamma}(J_{\beta}X,Y)+\sum_{i=1}^{4n}\epsilon_{i}\epsilon_{\gamma}R(Y,e_i,J_{\beta}X,J_{\gamma}e_i)+\sum_{i=1}^{4n}\epsilon_{i}\epsilon_{\gamma}R(e_i,J_{\beta}X,Y,J_{\gamma}e_i)=-2\epsilon_{\gamma}(\nabla_{J_{\beta}X}t)J_{\gamma}Y-\\
2\epsilon_{\beta}(\nabla_{Y}t)J_{\alpha}X+\sum_{i=1}^{4n}\epsilon_{i}\epsilon_{\gamma}(\nabla_{e_{i}}T)(J_{\beta}X,Y,J_{\gamma}e_i)+\sum_{i=1}^{4n}\epsilon_{i}\epsilon_{\gamma}g(T(J_{\beta}X,Y),T(e_i,J_{\gamma}e_i))-2\epsilon_{\gamma}L_{\gamma}(J_{\beta}X,Y).
\end{gather*}
\begin{gather*}
2\epsilon_{\gamma}\rho_{\gamma}(X,J_{\beta}Y)+\sum_{i=1}^{4n}\epsilon_{i}\epsilon_{\gamma}R(J_{\beta}Y,e_i,X,J_{\gamma}e_i)+\sum_{i=1}^{4n}\epsilon_{i}\epsilon_{\gamma}R(e_i,X,J_{\beta}Y,J_{\gamma}e_i)=2\epsilon_{\beta}(\nabla_{X}t)J_{\alpha}Y+\\
2\epsilon_{\gamma}(\nabla_{J_{\beta}Y}t)J_{\gamma}X+\sum_{i=1}^{4n}\epsilon_{i}\epsilon_{\gamma}(\nabla_{e_{i}}T)(X,J_{\beta}Y,J_{\gamma}e_i)+\sum_{i=1}^{4n}\epsilon_{i}\epsilon_{\gamma}g(T(X,J_{\beta}Y),T(e_i,J_{\gamma}e_i))-2\epsilon_{\gamma}L_{\gamma}(X,J_{\beta}Y).
\end{gather*}
The sum of all these equalities, $(\ref{11})$ and the fact that
$T$ is $(1,2)+(2,1)$-form with respect to each $J_{\alpha}$, gives
$$\frac{2(n-1)}{n}\rho_{\alpha}(J_{\beta}X,J_{\beta}Y)+\frac{2(n-1)}{n}\epsilon_{\beta}\rho_{\alpha}(X,Y)+\frac{2(n-1)}{n}\epsilon_{\gamma}\rho_{\gamma}(J_{\beta}X,Y)+\frac{2(n-1)}{n}\epsilon_{\gamma}\rho_{\gamma}(X,J_{\beta}Y)=$$
$$=-2L_{\alpha}(J_{\beta}X,J_{\beta}Y)-\epsilon_{\beta}2L_{\alpha}(X,Y)-\epsilon_{\gamma}2L_{\gamma}(J_{\beta}X,Y)-\epsilon_{\gamma}2L_{\gamma}(X,J_{\beta}Y).$$
From Lemma$\ref{21.1}$ and fact that $(n>1)$, we have
$(\ref{21.2})$ \hfill {\bf Q.E.D.}

We easily derive from Theorem$\ref{thm10}$
\begin{co}
The (2,0)+(0,2)-parts of the 2-forms $(dT)_{\alpha}$,
$(dT)_{\beta}$ with respect to $J_{\gamma}$ coincide.
\end{co}
\begin{thm}\label{thm10.1}
On a 4n-dimensional $(n>1)$ PQKT-manifold the following formula
hold
\begin{equation}\label{22.2}
\epsilon_{\alpha}\rho_{\alpha}(X,J_{\alpha}Y)+\epsilon_{\alpha}\rho_{\alpha}(J_{\alpha}X,Y)=-\frac{n}{n+1}(dt(X,Y)+\epsilon_{\alpha}dt(J_{\alpha}X,J_{\alpha}Y))
\end{equation}
In particular, $\rho_{\alpha}$ is of type $(1,1)$ with respect to
$J_{\alpha}$,$\alpha=1,2,3$ if and only if $dt$ is of type $(1,1)$
with respect to each $J_{\alpha}$,$\alpha=1,2,3$.
\end{thm}
{\bf Proof.} From the first Bianchi identity, formulas $(\ref{3})$
and $(\ref{11})$ follow
\begin{equation}\label{22.21}
2(\epsilon_{\alpha}\rho_{\alpha}(X,J_{\alpha}Y)+\epsilon_{\alpha}\rho_{\alpha}(J_{\alpha}X,Y))-\epsilon_{\alpha}(Ric(J_{\alpha}X,J_{\alpha}Y)-Ric(J_{\alpha}Y,J_{\alpha}X))-((Ric(X,Y)-Ric(Y,X))+
\end{equation}
$$+\frac{1}{n}(\epsilon_{\beta}\rho_{\beta}(X,J_{\beta}Y)+\epsilon_{\beta}\rho_{\beta}(J_{\beta}X,Y)+\epsilon_{\gamma}\rho_{\gamma}(X,J_{\gamma}Y)+\epsilon_{\gamma}\rho_{\gamma}(J_{\gamma}X,Y)-\rho_{\beta}(J_{\alpha}X,J_{\gamma}Y)-\rho_{\beta}(J_{\gamma}X,J_{\alpha}Y)+$$
$$+\rho_{\gamma}(J_{\alpha}X,J_{\beta}Y)+\rho_{\gamma}(J_{\beta}X,J_{\alpha}Y))=-2(dt(X,Y)+\epsilon_{\alpha}dt(J_{\alpha}X,J_{\alpha}Y))+\delta T(X,Y)+\epsilon_{\alpha}\delta T(J_{\alpha}X,J_{\alpha}Y).$$
First, we substitute $X \to J_{\alpha}X$ into $(\ref{21.2})$ to
\begin{equation}\label{22.22}
\epsilon_{\gamma}\rho_{\alpha}(J_{\gamma}X,J_{\beta}Y)+\epsilon_{\beta}\rho_{\alpha}(J_{\alpha}X,Y)+\rho_{\gamma}(J_{\gamma}X,Y)+\epsilon_{\gamma}\rho_{\gamma}(J_{\alpha}X,J_{\beta}Y)=0
\end{equation}
After that we substitute $Y \to J_{\alpha}Y$ into $(\ref{21.2})$
to get
\begin{equation}\label{22.23}
\epsilon_{\gamma}\rho_{\alpha}(J_{\beta}X,J_{\gamma}Y)+\epsilon_{\beta}\rho_{\alpha}(X,J_{\alpha}Y)+\rho_{\gamma}(X,J_{\gamma}Y)+\epsilon_{\gamma}\rho_{\gamma}(J_{\beta}X,J_{\alpha}Y)=0
\end{equation}
Summing up $(\ref{22.22})$ and $(\ref{22.23})$, we obtain
\begin{equation}\label{22.24}
\epsilon_{\alpha}\rho_{\alpha}(X,J_{\alpha}Y)+\epsilon_{\alpha}\rho_{\alpha}(J_{\alpha}X,Y)=\epsilon_{\gamma}\rho_{\gamma}(J_{\gamma}X,Y)+\epsilon_{\gamma}\rho_{\gamma}(X,J_{\gamma}Y)+
\end{equation}
$$+\rho_{\gamma}(J_{\alpha}X,J_{\beta}Y)+\rho_{\gamma}(J_{\beta}X,J_{\alpha}Y)+(\rho_{\alpha}(J_{\gamma}X,J_{\beta}Y)+\rho_{\alpha}(J_{\beta}X,J_{\gamma}Y))$$
We make the cyclic permutation $(\alpha, \beta, \gamma)$ $\to$
$(\beta, \gamma, \alpha)$ in $(\ref{22.24})$ to obtain
\begin{equation}\label{22.25}
\epsilon_{\alpha}\rho_{\alpha}(X,J_{\alpha}Y)+\epsilon_{\alpha}\rho_{\alpha}(J_{\alpha}X,Y)=\epsilon_{\beta}\rho_{\beta}(J_{\beta}X,Y)+\epsilon_{\beta}\rho_{\beta}(X,J_{\beta}Y)+
\end{equation}
$$-\rho_{\beta}(J_{\alpha}X,J_{\gamma}Y)-\rho_{\beta}(J_{\gamma}X,J_{\alpha}Y)-(\rho_{\alpha}(J_{\gamma}X,J_{\beta}Y)+\rho_{\alpha}(J_{\beta}X,J_{\gamma}Y))$$
Adding $(\ref{22.24})$ to $(\ref{22.25})$, we get
\begin{equation}\label{22.26}
2(\epsilon_{\alpha}\rho_{\alpha}(X,J_{\alpha}Y)+\epsilon_{\alpha}\rho_{\alpha}(J_{\alpha}X,Y))=\epsilon_{\gamma}\rho_{\gamma}(J_{\gamma}X,Y)+\epsilon_{\gamma}\rho_{\gamma}(X,J_{\gamma}Y)+
\end{equation}
$$+\epsilon_{\beta}\rho_{\beta}(J_{\beta}X,Y)+\epsilon_{\beta}\rho_{\beta}(X,J_{\beta}Y)+\rho_{\gamma}(J_{\alpha}X,J_{\beta}Y)+\rho_{\gamma}(J_{\beta}X,J_{\alpha}Y)-$$
$$-\rho_{\beta}(J_{\alpha}X,J_{\gamma}Y)-\rho_{\beta}(J_{\gamma}X,J_{\alpha}Y)$$
Now, equalities $(\ref{22.21})$, $(\ref{22.25})$, $(\ref{22.26})$
and $Ric(X,Y)-Ric(Y,X)=-\delta T(X,Y)$ (see \cite{I1}) prove the
assertion. \hfill {\bf Q.E.D.}
\begin{co}\label{10.2}
On a 4n-dimensional $(n>1)$ PQKT-manifold the following formula
hold
\begin{equation}\label{22.3}
\epsilon_{\alpha}\rho^g_{\alpha}(X,J_{\alpha}Y)+\epsilon_{\alpha}\rho^g_{\alpha}(J_{\alpha}X,Y)=-\frac{n-1}{2(n+1)}(dt(X,Y)+\epsilon_{\alpha}dt(J_{\alpha}X,J_{\alpha}Y))
\end{equation}
In particular, $\rho^{\star}_{\alpha}$ is symmetric if and only if
$dt$ is of type $(1,1)$ with respect to each
$J_{\alpha}$,$\alpha=1,2,3$.
\end{co}
{\bf Proof.} We get from $(\ref{15})$ that
$$\epsilon_{\alpha}\rho^g_{\alpha}(X,J_{\alpha}Y)+\epsilon_{\alpha}\rho^g_{\alpha}(J_{\alpha}X,Y)=\epsilon_{\alpha}\rho_{\alpha}(X,J_{\alpha}Y)+\epsilon_{\alpha}\rho_{\alpha}(J_{\alpha}X,Y)+\frac{1}{2}(dt(X,Y)+\epsilon_{\alpha}dt(J_{\alpha}X,J_{\alpha}Y)).$$
Now $(\ref{22.3})$ is a consequence of  $(\ref{22.2})$. \hfill
{\bf Q.E.D.}
\begin{pro}\label{orp1}
On a 4n-dimensional $(n>1)$ PQKT-manifold we have the equalities:
\begin{equation}\label{22.4}
Scal_{\alpha,\alpha}=Scal_{\beta,\beta}=Scal_{\gamma,\gamma},\quad
Scal_{\alpha,\beta}=0, \quad
Scal_{\alpha}=\frac{1}{2}(dt,\Phi_{\alpha})
\end{equation}
\end{pro}
{\bf Proof.}  Using $(\ref{22})$, we obtain
\begin{equation}\label{nov2}
\frac{2(n-1)}{n}(\rho_{\alpha}(X,J_{\alpha}Y)-
\rho_{\beta}(X,J_{\beta}Y))
=((dT)_{\alpha}(X,J_{\alpha}Y)-(dT)_{\beta}(X,J_{\beta}Y));
\end{equation}
\begin{equation}\label{22}
(n-1)\epsilon_{\alpha}\rho_{\alpha}(X,J_{\alpha}Y)=\frac{n(n-1)}{n+2}Ric(X,Y)-\frac{n(n-1)}{n+2}(\nabla_Xt)Y+
\end{equation}
$$+ \frac{n}{4(n+2)}\left\{(n+1)\epsilon_{\alpha}(dT)_{\alpha}(X,J_{\alpha}Y)-
\epsilon_{\beta}(dT)_{\beta}(X,J_{\beta}Y)
-\epsilon_{\gamma}(dT)_{\gamma}(X,J_{\gamma}Y)\right\}.\nonumber$$
Take the appropriate trace in (\ref{nov2}), to get
$Scal_{\alpha,\alpha}=Scal_{\beta,\beta},\quad
Scal_{\alpha,\beta}=0$.  The last equality in $(\ref{22.4})$ is a
direct consequence of $Scal_{\alpha,\beta}=0$ and $(\ref{22})$.
\hfill {\bf Q.E.D.}
\medskip

\noindent {\bf Definition.} {\em The three coinciding traces of
the Ricci forms on a $4n$ dimensional PQKT manifold $(n>1)$, give
a well-defined global function. We call this function the {\bf
paraquaternionic scalar curvature of the PQKT connection} and
denote it by $Scal_{\mathbb P}:= Scal_{\alpha,\alpha}$.}

\medskip

\begin{pro}\label{orp2}
On a 4n-dimensional $(n>1)$ PQKT manifold we have
\begin{equation}\label{pq1}
Scal^g_{\alpha}=Scal^g_{\beta}=Scal^g_{\gamma}= Scal_{\mathbb P} -
\delta t+||t||^2-\frac{1}{12}||T||^2,
-\epsilon_{\gamma}Scal^g_{\alpha,\beta}=Scal_{\gamma}=\frac{1}{2}(dt,\Phi_{\gamma}).
\end{equation}
\end{pro}
{\bf Proof.} The curvature $R^g$ of the
 Levi-Civita connection is related to $R$ via (\ref{15})
Taking the  traces in (\ref{15}) and using $(\ref{n1})$, we obtain
\begin{eqnarray}\label{rn2}
\epsilon_{\alpha}\rho^g_{\alpha}(X,J_{\alpha}Y)&=&\epsilon_{\alpha}\rho_{\alpha}(X,J_{\alpha}Y)+
\frac{1}{2}(\nabla_Xt)Y
-\epsilon_{\alpha}\frac{1}{2}(\nabla_{J_{\alpha}Y}t)J_{\alpha}X\\
\nonumber
&+&\frac{1}{2}\epsilon_{\alpha}t(J_{\alpha}T(X,J_{\alpha}Y))
+\frac{1}{4}\sum _{i=1}^{4n}
\epsilon_{i}\epsilon_{\alpha}g\left(T(X,e_i),T(J_{\alpha}Y,J_{\alpha}e_i)\right),
\end{eqnarray}
To finish, take the appropriate traces in (\ref{rn2}) and apply
Proposition~\ref{orp1}. \ \hfill {\bf Q.E.D.}

\medskip

\noindent {\bf Definition.}  {\em The three coinciding traces of
the Riemannian Ricci forms on a $4n$ dimensional PQKT manifold
$(n>1)$, give a well-defined global function. We call this
function the {\bf paraquaternionic $*$-scalar curvature} and
denote it by $Scal^g_{\mathbb P} := Scal^g_{\alpha,}$.}

\

\medskip

\begin{pro}\label{orp4}
On a 4n-dimensional ($n>1$) PQKT manifold $(M,g,\mathbb P)$ the
scalar curvatures are related by
\begin{eqnarray}
\nonumber Scal^g &=& \frac{n+2}{n}Scal_{\mathbb P} -3\delta t
+2||t||^2 -\frac{1}{12}||T||^2,\\ \nonumber Scal^g_{\mathbb P} &=&
Scal_Q - \delta t +||t||^2 -\frac{1}{12}||T||^2,\\ \nonumber Scal
&=& \frac{n+2}{n}Scal_{\mathbb P} -3\delta t +2||t||^2
-\frac{1}{3}||T||^2.\nonumber
\end{eqnarray}
\end{pro}
{\bf Proof.} We derive from (\ref{15}) that
\begin{eqnarray}\label{r5}
Ric^g(X,Y) &=& Ric(X,Y) + \frac{1}{2}\delta T(X,Y) +
\frac{1}{4}\sum _{i=1}^{2n} g\left(T(X,e_i),T(Y,e_i)\right), \\ &
& Scal^g=Scal + \frac{1}{4}||T||^2.\nonumber
\end{eqnarray}
Take the trace in (\ref{22}) to get the first equality of the
proposition. The second equality is already proved in
Proposition~\ref{pq1}. The last one is a consequence of (\ref{r5})
and the already proven first equality in the proposition. \hfill
{\bf Q.E.D.}
\\
\section{PQKT manifolds with parallel torsion and homogeneous PQKT structures}

Let $(G/K,g)$ be a reductive (locally) homogeneous
pseudo-Riemannian manifold. The canonical connection $\nabla$ is
characterized by the properties $\nabla g=\nabla T=\nabla R=0$
\cite{K2}. A homogeneous paraquaternionic Hermitian manifold
(resp. homogeneous hyper para-Hermitian) manifold $(G/K,g,\mathbb
P)$ is a homogeneous pseudo-Riemannian manifold with an invariant
paraquaternionic Hermitian subbundle $\mathbb P$ (resp.  three
invariant anti commuting (para)complex structures ). This means
that the bundle $\mathbb P$ (resp. each of the (para)complex
structures) is parallel with respect to the canonical connection
$\nabla$. The torsion of $\nabla$ is totally skew-symmetric if and
only if the homogeneous pseudo-Riemannian manifold is naturally
reductive. Homogeneous PQKT (resp. HPKT) manifolds are homogeneous
paraquaternionic Hermitian (resp. homogeneous hyper
para-Hermitian) manifold which are naturally reductive.
\\
We show that there are no homogeneous PQKT manifold with torsion
4-form $dT$ of type (2,2) with respect to each $J_{\alpha}$ in
dimensions greater than four. First, we prove the following
technical result
\begin{pro}\label{teh}
Let $(M,g,(J_{\alpha}),\nabla)$ be a 4n-dimensional $(n>1)$ PQKT
manifold with 4-form $dT$ of type (2,2) with respect to each
$J_{\alpha}, \alpha =1,2,3$. Suppose that the torsion is parallel
with respect to the PQKT-connection. Then the Ricci forms
$\rho_{\alpha}$
 are given by
\begin{equation}\label{27}
\epsilon_{\alpha} \rho_{\alpha}(X,J_{\alpha}Y) = \lambda g(X,Y),
\quad \alpha =1,2,3,
\end{equation}
where $\lambda$ is a smooth function on $M$.
\end{pro}
{\it Proof.} Let the torsion be parallel i.e. $\nabla T=0$. This
implies that the Ricci tensor is symmetric \cite{F2}. The
equalities (\ref{13}) and (\ref{14}) lead to
\begin{equation}\label{17}
B(X,Y,Z,U)=  {\sigma \atop XYZ}
\left\{g(T(X,Y),T(Z,U)\right\}=\frac{1}{2}dT(X,Y,Z,U).
\end{equation}
We get $D=0$ from (\ref{tir1}).

Suppose now that the 4-form $dT$ is of type (2,2) with respect to
each $J_{\alpha}, \alpha =1,2,3.$. Then it satisfies the
equalities
\begin{equation}\label{26.9}
- \epsilon_{\alpha} dT(X,Y,Z,U)=dT(J_{\alpha}X,J_{\alpha}Y,Z,U) +
dT(J_{\alpha}X,Y,J_{\alpha}Z,U) + dT(X,J_{\alpha}Y,J_{\alpha}Z,U).
\end{equation}.
\\
The similar arguments as we used in the proof of
Proposition~\ref{l1} but
 applying
(\ref{26.9}) instead of (\ref{3}), yield
\begin{lem}\label{l2}
On a PQKT manifold with 4-form $dT$ of type (2,2) with respect to
each $J_{\alpha}, \alpha=1,2,3$, the following equalities hold:
\begin{equation}\label{24}
(dT)_1(X,J_1Y) = (dT)_2(X,J_2Y) = -(dT)_3(X,J_3Y),
\end{equation}
\begin{equation}\label{24'}
(dT)_{\alpha}(X,J_{\alpha}Y)= -(dT)_{\alpha}(J_{\alpha}X,Y), \quad
\alpha=1,2,3.
\end{equation}
\end{lem}
We substitute (\ref{24}), (\ref{17}) and $D=0$ into (\ref{22}) to
get
\begin{equation}\label{a1}
\rho_1(X,J_1Y) = \rho_2(X,J_2Y) =- \rho_3(X,J_3Y),
\end{equation}
\begin{equation}\label{25}
\epsilon_{\alpha}\rho_{\alpha}(X,J_{\alpha}Y)=
\frac{n}{n+2}Ric(X,Y) + \frac{n}{4(n+2)}\epsilon_{\alpha}
(dT)_{\alpha}(X,J_{\alpha}Y), \quad \alpha=1,2,3.
\end{equation}
The equality (\ref{24'}) shows that the 2-form $dT_{\alpha}$ is a
(1,1)-form with respect to $J_{\alpha}$. Hence, the $dT_{\alpha}$
is (1,1)-form with respect to each $J_{\alpha},\alpha=1,2,3$,
because of (\ref{24}). Since the Ricci tensor $Ric$ is symmetric,
(\ref{25}) shows that the Ricci tensor  $Ric$ satisfies
$Ric(J_{\alpha}X,J_{\alpha}Y)=-\epsilon_{\alpha} Ric(X,Y), \alpha
=1,2,3$ for each $J_{\alpha}$ and  the Ricci forms $\rho_{\alpha},
\alpha =1,2,3$ are (1,1)-forms with respect to all $J_{\alpha},
\alpha=1,2,3$. Taking into account (\ref{11}), we obtain
\begin{equation}\label{26}
-\epsilon_{\alpha}R(X,J_{\alpha}X,Z,J_{\alpha}Z)+R(X,J_{\alpha}X,J_{\beta}Z,J_{\gamma}Z)+
R(J_{\beta}X,J_{\gamma}X,Z,J_{\alpha}Z)-
\end{equation}
$$-\epsilon_{\alpha}R(J_{\beta}X,J_{\gamma}X,J_{\beta}Z,J_{\gamma}Z)=\frac{1}{n}\left(-\epsilon_{\alpha}\rho_{\alpha}(X,J_{\alpha}X)+
\rho_{\alpha}(J_{\beta}X,J_{\gamma}X)\right) g(Z,Z) =$$
$$=-\frac{2}{n}\epsilon_{\alpha} \rho_{\alpha}(X,J_{\alpha}X)g(Z,Z)$$
where the last equality of (\ref{26}) is a consequence of the
following identity
$$
\rho_{\alpha}(J_{\beta}X,J_{\gamma}X)=\epsilon_{\beta}\rho_{\beta}(J_{\beta}X,X)=
-\epsilon_{\alpha}\rho_{\alpha}(X,J_{\alpha}X).
$$
The left hand side of (\ref{26}) is symmetric with respect to the
vectors $X,Z$ because $D=0$. Hence,
$\rho_{\alpha}(X,J_{\alpha}X)g(Z,Z)=\rho_{\alpha}(Z,J_{\alpha}Z)g(X,X),
\alpha=1,2,3$. The last equality together with (\ref{a1}) implies
(\ref{27}). \hfill {\bf Q.E.D.}
\begin{thm}\label{thm2}
Let $(M,g,(J_{\alpha}))$ be a 4n-dimensional ($n>1$) PQKT manifold
with 4-form $dT$ of type (2,2) with respect to each $J_{\alpha},
\alpha =1,2,3$. Suppose that the torsion is parallel with respect
to the PQKT-connection. Then  $(M,g,(J_{\alpha}))$ is either a
HPKT manifold with parallel torsion or a PQK manifold.
\end{thm}
{\it Proof.} We apply Proposition~\ref{teh}. If the function
$\lambda =0$ then $\rho_{\alpha}=0, \alpha =1,2,3$, by (\ref{27})
and Proposition~\ref{p2} implies that the PQKT manifold is
actually a HPKT manifold.

 Let $\lambda \not= 0$. The condition (\ref{27}) determines the torsion
completely. We proceed involving (\ref{12}) into the computations.
We calculate, using ({\ref{1}) and (\ref{27}), that
\begin{equation}\label{28}
(\nabla_Z\rho_{\alpha})(X,Y)
=\lambda\left\{\omega_{\beta}(Z)F_{\gamma}(X,Y) +
\epsilon_{\gamma}\omega_{\gamma}(Z)F_{\beta}(X,Y)\right\} +
d\lambda(Z)F_{\alpha}(X,Y).
\end{equation}
Applying the operator $d$ to (\ref{11}), we get taking into
account (\ref{27})
 that
\begin{equation}\label{29}
d\rho_{\alpha}(X,Y,Z) = \lambda(\epsilon_{\gamma} F_{\beta}\wedge
\omega_{\gamma} + \omega_{\beta} \wedge F_{\gamma})
\end {equation}
On the other hand, we have
\begin{equation}\label{30}
d\rho_{\alpha} = {\sigma \atop XYZ}
\left\{(\nabla_Z\rho_{\alpha})(X,Y)
+\lambda(T(X,Y,J_{\alpha}Z)\right\}, \quad \alpha=1,2,3.
\end{equation}
Comparing the left-hand sides of (\ref{29}) and (\ref{30}) and
using (\ref{28}), we derive $$ \lambda {\sigma \atop
XYZ}\left\{g(T(X,Y),J_{\alpha}Z)\right\} = - d\lambda\wedge
F_{\alpha}(X,Y,Z), \quad \alpha=1,2,3. $$ The last equality
implies $ \lambda T=- \epsilon_{\alpha} J_{\alpha}d\lambda \wedge
F_{\alpha}, \quad \alpha=1,2,3. $ If $\lambda$ is a non zero
constant then $T=0$. If $\lambda$ is not a constant then there
exists a point $p\in M$ and a neighbourhood $V_p$ of $p$ such that
$\lambda \Big |_{V_p} \not=0$. Then
\begin{equation}\label{c10}
T= - \epsilon_{\alpha} J_{\alpha}d\ln\lambda \wedge F_{\alpha},
\quad \alpha=1,2,3.
\end{equation}
We take  the trace in (\ref{c10}) to obtain
\begin{equation}\label{c11}
4(n-1)J_{\alpha}d\ln\lambda = 0, \quad \alpha =1,2,3.
\end{equation}
The equation (\ref{c11}) forces $d\lambda =0$ since $n>1$ and
consequently $T=0$ by (\ref{c10}).  Hence, the QKT space is a QK
manifold which completes the proof. \hfill {\bf Q.E.D.}

On a locally homogeneous PQKT manifold the torsion and curvature
are parallel and Theorem~\ref{thm2} leads to the following
\begin{thm}\label{q1}
A (locally) homogeneous 4n-dimensional $(n>1)$ PQKT manifold with
torsion 4-form $dT$ of type (2,2) is either (locally) homogeneous
HPKT space or a (locally) symmetric PQK space.
\end{thm}
Theorem~\ref{q1} shows that there are no homogeneous (proper) PQKT
manifolds with torsion 4-form of type (2,2) in dimensions greater
than four.

\bibliographystyle{hamsplain}

\providecommand{\bysame}{\leavevmode\hbox
to3em{\hrulefill}\thinspace}






\end{document}